# Variational problems for Riemannian functionals and arithmetic groups.

## Alexander Nabutovsky and Shmuel Weinberger

In this paper we introduce a new approach to variational problems on the space $Riem(M^n)$ of Riemannian structures (i.e. isometry classes of Riemannan metrics) on any fixed compact manifold $M^n$ of dimension $n \geq 5$. This approach often enables one to replace the considered variational problem on $Riem(M^n)$ (or on some subset of $Riem(M^n)$) by the same problem but on spaces $Riem(N^n)$ for every manifold $N^n$ from a class of compact manifolds of the same dimension and with the same homology as $M^n$ but with the following two useful properties: (1) If $\nu$ is any Riemannian structure on any manifold $N^n$ from this class such that $Ric_{(N^n,\nu)} \geq -(n-1)$, then the volume of $(N^n, \nu)$ is greater than one; and (2) Manifolds from this class do not admit Riemannian metrics of non-negative scalar curvature. The first property is obviously helpful when one knows how to prove uniform diameter and curvature bounds for a minimizing sequence but wants to ensure that this sequence does not collapse to a metric space of lower dimension. The second property seems likely to be useful when one studies a variational problem on the space of Riemannian structures of constant scalar curvature (compare [An1], [An 2], [Sch]).

As a first application we prove a theorem which can be informally explained as follows: Let $M$ be any compact connected smooth manifold of dimension greater than four, $Met(M)$ be the space of isometry classes of compact metric spaces homeomorphic to $M$ endowed with the Gromov-Hausdorff topology, $Riem_1(M) \subset Met(M)$ be the space of Riemannian structures on $M$ such that the absolute values of sectional curvature do not exceed one, and $R_1(M)$ denote the closure of $Riem_1(M)$ in $Met(M)$. Then diameter regarded as a functional on $R_1(M)$ has infinitely many "very deep" local minima. Moreover, the set of its values at these local minima is unbounded. (This result can be regarded as a solution in dimensions $\geq 5$ of the Problem 1 in S.T. Yau's list of problems in differential geometry ([Y]) which asks for a way to define canonical Riemannian metrics on all compact manifolds (of a fixed dimension). A detailed discussion of this problem can be found in [Br].) We also give an exponentially growing lower bound for the distribution function of these deep local minima. These results are corollaries of our results about sublevel sets of the diameter regarded as a functional on $R_1(M)$. In particular, we demonstrate that for all sufficiently large $x$ the sublevel sets of the diameter $diam^{-1}((0,x]) \subset R_1(M)$ are not connected; these sublevel sets can be represented as a union of at least $\exp(c(n)x^n)$ non-empty subsets separated from each other by "gaps", where $c(n) > 0$ depends only on $n$; the infimum of the volume on some of these subsets is positive, and the number of the subsets with this property also grows at least exponentially with $x^n$; and this "severe" disconnectedness cannot be avoided by allowing a "controllable" increase of diameter along the path in $R_1(M)$ which connects Riemannian metrics from different connected components of $diam^{-1}((0,x])$.

We then try to apply these ideas to the question of existence of Einstein metrics (or

The first author was partially supported by the NSERC grant OGP 0155879; the second author was partially supported by an NSF grant.

at least almost Einstein metrics) of non-positive scalar curvature on $S^n$. We propose an approach which eventually might lead to the construction of such metrics. In the meantime we use our technique to demonstrate that the (contractible; cf. [L]) space of Riemannian metrics of constant negative scalar curvature on $S^n$, $n \geq 5$, has very complicated geometry (see Theorem 2 in section 1).

All these results are based on the existence of $n$-dimensional smooth homology spheres, ($n \geq 5$), not admitting Riemannian metrics of non-negative scalar curvature and such that the volume of any of these homology spheres with respect to any Riemannian metric of Ricci curvature $\geq -(n-1)$ is greater than one. The fundamental groups of these homology spheres in our construction are "made" of certain arithmetic groups with appropriate homology properties.

## 0. Introduction.

Let $M$ be a compact manifold. Assume that we are trying to prove the existence of a solution of some variational problem on the space of Riemannian structures on $M$ (or on some subset of this space). Choose some minimizing sequence $\{\mu_i\}_{i=1}^{\infty}$. According to the Gromov-Cheeger compactness theorem if we are able to prove the existence of uniform upper bounds for the absolute values of the sectional curvature and diameters as well as a uniform positive lower bound for the volumes of $(M, \mu_i)$, then we will be guaranteed that a subsequence of the sequence $\{\mu_i\}$ converges to a $C^{1,\alpha}$-smooth limit, where $\alpha$ is any positive number less than one. Even if we are able to prove only the existence of a uniform lower bound for the values of the sectional curvature (as well as the bounds for diameter and volume) instead of the upper bound for the absolute values of the sectional curvature, a compactness theorem proven by Burago, Gromov and Perelman ([BGP]) guarantees that a subsequence of $\{(M, \mu_i)\}$ converges in the Gromov-Hausdorff topology to an Alexandrov space. Since the interesting variational problems on the space of Riemannian structures are usually scale invariant, we can have an upper bound for the diameter for free just by an appropriate rescaling. To find a uniform two-sided or even lower bound for the sectional curvature is a non-trivial (and sometimes very difficult) problem which can be compared with deriving a priori bounds while trying to prove the existence of solutions of a PDE. However, even if we have this two-sided bound for the sectional curvature, the sequence $\{(M, \mu_i)\}$ can still collapse to some metric space of a lower dimension. One of the main conclusions of our paper is that if one is interested in *local* extrema then in a wide class of such situations the trouble caused by the absence of a positive uniform lower bound for the volume can be avoided: We demonstrate how a search for local minima of a functional on the space of Riemannian structures on an arbitrary Riemannian manifold $M_0^n$ can be often replaced by a search for local minima of the same functional but on spaces of Riemannian structures on every manifold $M^n$ such that $H^*(M^n) = H^*(M_0^n)$ but such that $vol_\mu(M^n) \geq w(n)$ for any Riemannian metric $\mu$ on $M^n$ satisfying $Ric_{(M^n, \mu)} \geq -(n-1)$, where $w(n)$ is a positive constant depending only on $n$. (In fact, one can even demand here that $vol_\mu(M^n) \geq C$, where $C$ is an arbitrary positive constant, e.g. $C = 1$.) Our method enables one to restrict even further the class of manifolds which one needs to consider in order to prove the existence of local extrema for $M$ even further. For example, it is sufficient to consider only manifolds which do not admit Riemannian metrics of non-negative scalar curvature (this property can be useful when one is looking for Einstein metrics). Moreover, when our method works, it yields an infinite set of distinct local minima and a (usually exponential) lower bound for the distribution function of these local minima. As a concrete application of our technique we prove that for any compact manifold $M$ of dimension $\geq 5$ the diameter has infinitely many local minima on the closure with respect to the Gromov-Hausdorff topology of the set of all smooth Riemannian structures on $M$ satisfying $|K| \leq 1$ in the space $Met(M)$ of all isometry classes of metric spaces homeomorphic to $M$ (see Theorem B below). ($Met(M)$ is assumed to be endowed with the Gromov-Hausdorff topology. It is known that all elements of this closure are $C^{1,\alpha}$-smooth Riemannian structures on $M$ for any $\alpha \in (0,1)$; cf. [F]). Note also that according to Theorem C at the end of section 1 the method developed in the present paper can be applied in a similar fashion not only to extremal problems but to a very wide class of (finite) systems of equations and inequalities for functionals on the space of Riemannian structures on any compact manifold of dimension $\geq 5$. Because our method involves results from several branches of mathematics we will be somewhat lengthy in our introduction.

To state our main theorems let us introduce the following notations: Let $M$ be any compact manifold. Denote by $Riem_1(M)$ the space of Riemannian structures on $M$ such that the absolute values of sectional curvature do not exceed one. (The basic properties of the space of Riemannian structures on a compact manifold can be found, for example, in [Bou] and [Be].) Denote by $Met(M)$ the space of isometry classes of metric spaces homeomorphic to $M$ endowed by the Gromov-Hausdorff metric, and by $R_1(M)$ the closure of $Riem_1(M)$ in $Met(M)$. Note that for many manifolds $M$ (including tori of any dimension and spheres of any odd dimension) $R_1(M)$ is not complete because of the possibility of a collapsing to a lower-dimensional metric space. On the other hand the Cheeger-Gromov compactness theorem implies that all elements of $R_1(M)$ are $C^{1,\alpha}$-smooth Riemannian structures on $M$, for any positive $\alpha < 1$. We can consider diameter as a functional on $R_1(M)$ which will be denoted by $diam$. An increasing function $\phi : (0, \infty) \longrightarrow (0, \infty)$ is said to be *effectively majorizable* if there exists a (Turing) computable increasing function $\alpha : \mathbb{N} \longrightarrow \mathbb{N}$ such that for any $x$ $\phi(x) \leq \alpha([x])$. (For example, for any positive real $c_1, c_2$ functions $c_1 x^{c_2}$, $\exp(c_1 \, x^{c_2})$ or $\exp(\exp(\ldots \exp(x^{c_1})))$ ($[x^{c_2}] + 1$ exponentiations) are effectively majorizable.

**Theorem A.(Informal version)** Let $M^n$ be any compact $n$-dimensional manifold of dimension $n \geq 5$. For all sufficiently large $x$ the sublevel set $diam^{-1}((0,x])$ of diameter regarded as a functional on $R_1(M^n)$ is not connected. Moreover, it can be represented as a union of disjoint non-empty subsets ("components") separated by a "gap" (in the Gromov-Hausdorff metric) and such that any two Riemannian structures from different components cannot be connected by a finite sequence of sufficiently small "jumps" in $R_1(M^n)$ passing through Riemannian structure of diameter $\leq x$ or even $\leq \phi(x)$, where $\phi$ is any increasing effectively majorizable function. (But the value of $x$ such that the last statement becomes true starting from this value depends on the choice $\phi$, of course.) The number of these "components" grows at least exponentially with $x^n$. Finally, (and most important for applications we have in mind!) for all sufficiently large $x$ for some of these "components" the volume regarded as a functional on a fixed "component" has a positive lower bound not less than a certain positive constant depending only on $n$. The number of "components" with this property also grows at least exponentially with $x^n$.

**Theorem A.(Formal version)** Let $M^n$ be any compact $n$-dimensional manifold, $n \geq 5$. Let $\phi$ be any effectively majorizable function such that $\phi(x) \geq x$ for any positive $x$. For any $x$ there exist $I \geq 1$, $J \in \{0, \ldots, I\}$, and a partition of the sublevel set $diam^{-1}((0, x]) \subset R_1(M^n)$ into disjoint non-empty subsets $R_i(M^n, x)$, $i = 1, \ldots, I$ with the following properties:
(i) There exists a positive $E(n)$ depending only on $n$ such that the Gromov-Hausdorff distance between any Riemannian metric from $R_i(M^n, x)$ and $R_j(M^n, x)$ for $i \neq j$ is not less than $E(n) \exp(-(n-1)x)$;
(ii) Moreover, for any $i, j \in \{1, \ldots, I\}$, $i \neq j$ and any $\mu \in R_i(M^n, x), \nu \in R_j(M^n, x)$ there is no sequence of "jumps" in $R_1(M^n)$ of length not exceeding $E(n) \exp(-(n-1)\phi(x))$ connecting $\mu$ and $\nu$ and passing only through Riemannian structures on $M^n$ from $diam^{-1}((0, \phi(x)])$. (That is, there is no finite sequence of elements of $R_1(M^n)$ of diameter not exceeding $\phi(x)$ such that the first term of the sequence coincides with $\mu$, the last term of the sequence coincides with $\nu$ and the Gromov-Hausdorff distance between two consecutive terms of the sequence does not exceed $E(n) \exp(-(n-1)\phi(x))$.);
(iii) There exist a positive constant $v(n)$ depending only on $n$ such that for any positive integer $j \leq J$ and for any $\mu \in R_j(M^n, x)$ the volume of $(M^n, \mu)$ is not less than $v(n)$;
(iv) For a strictly increasing unbounded sequence $\{x_i\}_{i=1}^\infty$ of values of $x$ $I \geq 2$ and $J \geq 1$. (This implies, in particular, that for any $i$, $diam^{-1}((0, x_i]) \subset R_1(M^n)$ is not connected and that $\inf_{\mu \in R_1(M^n, x_i)} vol_\mu(M^n) \geq v(n)$.);
(v) Furthermore, $J \geq 1$ and $I \geq 2$ for all sufficiently large $x$ (and not only for some infinite unbounded sequence $\{x_i\}_{i=1}^\infty$ of values of $x$). Moreover, there exists a constant $C(n) > 0$ depending only on $n$ such that for all sufficiently large $x$ $I \geq J \geq \exp(C(n)x^n)$.

**Remark 1.** It is clear that Theorem A without (iv) or (v) is vacuous. Of course, (v) is stronger than (iv), and the second statement in (v) is stronger than the first statement. The reason why we decided to state (iv) separately is that its proof is conceptually and technically simpler, and already has quite interesting applications (see Remark 3 after the text of Theorem B below). Moreover, combining the proof of parts (i)-(iv) of Theorem A given below with Lemma 6 in [N3] in a fashion similar to its application in [N3] one immediately obtains the first statement of (v) (that is, the statement that $I \geq 2$ and $J \geq 1$ for all sufficiently large $x$.) The proof of the exponential growth of $I$ and $J$ given below uses the notion of time-bounded Kolmogorov complexity and heavily relies on the constructions and results of [N2]. Also, observe that (i) is a particular case of (ii), where $\phi(x) = x$. Here the purpose of formulating first a weaker property (i) was to state a result which a reader not willing to deal with the terminology from recursion theory can use. (However, we do not know a proof of the existence of the partition satisfying (i),(iii),(iv) essentially simpler than the proof of the existence of a partition satisfying (i)-(iv) given below.)

**Important remark 2.** Theorem A can be strengthened in the following two ways:
(1) In (iii) one can replace the statement that the infimum of volume on $R_j(M^n, x)$ for any $j \in \{1, \ldots J\}$ is not less than the positive constant $v(n)$ depending only on $n$ by the statement that the infimum of volume on $R_j(M^n, x)$ is not less than $\exp(q(n)x)$, where $q(n)$ is a positive constant depending only on $n$. As a corollary, we can also replace $v(n)$ in the text of Theorem A by any positive constant of our choice, e.g. by 1. (As before, part (v) of Theorem A contains the statement that the number of $j$ with this property is at least $\exp(C(n)x^n)$ for all sufficiently large $x$.)

(2) In Theorem A, part (iii) one can replace the condition that for any $\mu \in R_j(M^n, x)$, $j \in \{1, \ldots J\}$, the volume of $(M^n, \mu)$ is not less than $v(n)$ by the condition that the the volume of the subset of $(M^n, \mu)$ which consists of all points where the injectivity radius is greater than some positive constant $\epsilon(n)$ depending only on $n$ is not less than $v(n)$. Furthermore, combining both parts of this remark together we can replace $v(n)$ here by $\exp(q(n)x)$, where $q(n)$ is the same as above.

This remark will be proven at the end of section 4 after the text of the proof of Theorem A.

For any continuous functional the local minima of its restriction to its sublevel sets will automatically be its local minima on the whole space. Therefore local minima of diameter on its sublevel sets will be automatically its local minima on $R_1(M^n)$. The Gromov-Cheeger compactness theorem implies the existence of the global minimum of diameter on $R_j(M^n, x)$ for any $j \leq J$. These global minima are, of course, local minima of $diam$ on $diam^{-1}((0, x])$. Moreover, these local minima are $C^{1,\alpha}$-smooth Riemannian structures on $M^n$ for any $\alpha \in (0, 1)$. Thus, applying Theorem A we immediately obtain the following result:

**Theorem B.** For any compact manifold $M^n$ of dimension $n \geq 5$ the set of local minima of diameter regarded as a functional on $R_1(M^n)$ is infinite, and the set of values of diameter at its local minima is unbounded. Furthermore, let $\phi$ be any effectvely majorizable function such that $\phi(x) \geq x$ (e.g. $\phi(x) = x$). Then for any sufficiently large $x$ there exists at least $[\exp(C(n)x^n)]$ local minima of $diam$ on $R_1(M^n)$ such that the value of diameter at any of these local minima does not exceed $x$, the volume is not less than $v(n) > 0$, and which are "deep" in the following sense: Let $\mu$ be one of these local minima. There is no finite sequence of "jumps" of length $\leq E(n) \exp(-(n-1)\phi(diam(\mu)))$ in $R_1(M^n)$ connecting $\mu$ with either a Riemannian structure on $M^n$ of a smaller diameter or with another of these local minima and passing only through Riemannian structures of diameter $\leq \phi(diam(\mu))$ on $M^n$. (The constants $C(n)$ and $E(n)$ here are the same as in Theorem A.)

**Remark 1.** This theorem can be regarded as a solution in dimension greater than or equal to five of the following problem which appears as Problem 1 in the S.T. Yau list of problems in differential geometry [Y]: *Find a general way to construct canonical Riemannian metrics on a given compact manifold.* (See also [Br], [Sr].) The considered variational problem is quite natural (see Remark 4 below) and has "not very large" set of solutions for *every* compact $n$-dimensional manifold, ($n \geq 5$). Of course, it would be desirable to obtain a better smoothness of the canonical Riemannian metrics than just $C^{1,\alpha}$ for any $\alpha < 1$. On the other hand, the non-uniqueness apparently cannot be completely avoided: According to Theorem 3 of [N1] if one insists that a canonical Riemannian metric must exist on any smooth compact manifold of the considered dimension $n \geq 5$, and if there is an algorithm recognizing a Riemannian metric sufficiently close to a canonical as having this property, then the set of canonical Riemannian structures of volume one on any compact manifold of the considered dimension is infinite. The condition of the existence of the algorithm recognizing Riemannian metrics sufficiently close to canonical is quite natural here: Without this condition one can just choose one Riemannian metric on any compact manifold using the axiom of choice.

**Remark 2.** Note that for many manifolds the *global* minimum of diameter on $R_1(M)$ does not exist. For example, the infimum of diameter on $R_1(T^n)$ is equal to zero. (Here $T^n$ denotes the $n$-dimensional torus.)

**Remark 3.** The proof of the exponential lower bound for $I$, $J$ in part (v) of Theorem A heavily relies on the material of sections 2,3 of [N2]. For the reader not intrested in these exponential lower bounds note that parts (i)-(iv) of Theorem A imply the following weaker version of Theorem B: *Let $M^n$ be as before a compact connected manifold of dimension $n \geq 5$. Let $\phi$ be an effectively majorizable function such that $\phi(x) \geq x$. There exists an infinite sequence $\{\mu_i\}_{i=1}^{\infty}$ of local minima of diam on $R_1(M^n)$ such that the corresponding sequence of values of diameter $\{diam(\mu_i)\}_{i=1}^{\infty}$ is unbounded; and for any $i$ there is no sequence of "jumps" of length $\leq E(n)\exp(-(n-1)\phi(diam(\mu_i)))$ in $diam^{-1}((0, \phi(diam(\mu_i))]) \subset R_1(M^n)$ connecting $\mu_i$ with a Riemannian metric on $M^n$ of diameter strictly less than $diam(\mu_i)$.* Indeed, for an unbounded sequence $\{x_i\}_{i=1}^{\infty}$ of values of $x$ $J \geq 1$ and $I \geq 2$. We can proceed as follows. Start from $X_1 = x_1$. Define $\mu_1$ as the global minimum of $diam$ on $R_1(M^n, X_1)$. For any $i \geq 2$ choose $X_i \in \{x_i\}_{i=1}^{\infty}$ to be sufficiently large to ensure that $diam^{-1}((0, X_{i-1}+1])$ belongs to $R_{k(i)}(M^n, X_i)$ for some $k(i)$. (The existence of such $X_i$ follows from the precompactness of $diam^{-1}((0, X_{i-1}+1])$.) If $M^n$ admits Riemannian metrics such that $\sup|K| \leq 1$ and $vol < v(n)$, then we can choose $X_i$ to be sufficiently large to ensure that $\inf_{\mu \in R_{k(i)}(M^n, X_i)} vol_\mu(M^n) < v(n)$. (Here $v(n)$ is the same as in Theorem A, (iii).) Let $l = 1$ in this case. If $\inf_{\mu \in R_1(M^n)} vol_\mu(M^n) \geq v(n)$, then $J = I \geq 2$. Choose (any) $l \in \{1, 2\}$ such that $l \neq k(i)$. In both cases there exists the global minimum $\mu_i$ of $diam$ on $R_l(M^n, X_i)$, and the value of the diameter at this global minimum is greater than $X_{i-1} + 1$. It is easy to see that the sequence $\{\mu_i\}$ has the required properties.

**Remark 4.** If $M^n$ does not admit a flat metric then the local minima of *diam* on $R_1(M^n)$ are the same as local minima of $\sup|K|$ on the space of Riemannian structures of diameter one on $M^n$, where $\sup|K|$ is understood in the sense of Alexandrov spaces with curvature bounded from both sides (cf. [BN]). In this setting the variational problem considered in Theorem B looks very natural.

**Remark 5.** The supplements to Theorem A stated as Remark 2 after its statement imply that in the statement of Theorem B we can also demand that for any of the local minima $\mu$ $vol_\mu(M^n) \geq \exp(q(n)x)$. Moreover, we can demand that the volume of the subset of $(M^n, \mu)$ formed by all points where the injectivity radius is greater than $\epsilon(n) > 0$ is greater than $\exp(q(n)x)$. Here $q(n)$ and $\epsilon(n)$ are the same as in the text of Remark 2 after Theorem A.

Here is a brief description of our approach. Its central and the most difficult part is the construction for any $n \geq 5$ of a smooth homology sphere not admitting Riemannian metrics of non-negative scalar curvature and such that the volume of this homology sphere endowed with any Riemannian metric of Ricci curvature $\geq -(n-1)$ is not less than a certain positive constant $w(n)$ depending only on $n$. (These properties are typical for manifolds admitting non-flat Riemannian metrics of non-positive sectional curvature. However, it is not presently known whether or not there exist homology $n$-spheres admitting metrics of non-positive sectional curvature or even $K(\pi, 1)$ smooth homology $n$-spheres for $n > 3$.) We build the fundamental groups of these homology spheres out of certain arithmetic groups by taking the universal central extension of an appropriate amalgamated free product with several copies of an acyclic group. Our constructions heavily rely on deep results of several mathematicians on cohomology of arithmetic groups ([Bor 2], [Bor 3], [BorWal], [Cl0], [Cl], [Kt0], [Kt1]). Other crucial ingredients are the classification of smooth homology spheres by J.-P. Hausmann and P. Vogel (cf. [Haus], [V]) based on S. Cappell and J. Shaneson's theory of homology surgery ([CS]) and M. Gromov's theorems establishing for a class of manifolds that the volume of any of these manifolds endowed with any Riemannian metric of Ricci curvature $\geq -(n-1)$ is not less than a certain positive constant depending only on the dimension (section 6.6 in [Gr 2] and the Main Inequality on p.12 of [Gr 1]). The required homology spheres are constructed in section 2, modulo Proposition 4 establishing the existence of finitely presented groups with certain homology properties. In section 3 we describe several constructions of such groups.

Then we prove that for any fixed $n \geq 5$ there is no algorithm which distinguishes such smooth homology spheres from the standard sphere $S^n$. This result is stated as Theorem 1.A at the beginning of section 1 and proven in section 2 using Proposition 4. (Very roughly, the idea is that starting from one group which is the fundamental group of a smooth homology sphere with the required properties one can construct infinitely many of such groups. In fact, there is enough flexibility in construction of fundamental groups of such homology spheres to prove that there is no algorithm distinguishing such groups from the trivial group. The proof of this fact uses ideas of the proof of the classical Adyan-Rabin theorem establishing the algorithmic unsolvability of the triviality problem for finitely presented group given in [M]. By the way, observe that even if there exist hyperbolic smooth homology spheres of dimensions $\geq 5$ the problem of distinguishing between such spheres and the standard sphere is algorithmically solvable.) Furthermore, Theorem 1 provides also a lower bound for the time-bounded Kolmogorov complexity of the problem of recognition of $S^n$ among homology spheres from an appropriate effectively constructed sequence for any computable time bound. (Any homology sphere in this sequence is either diffeomorphic to $S^n$ or does not admit Riemannian metrics of non-negative scalar curvature and does not admit Riemannian metrics such that $Ric \geq -(n-1)$ and the volume $< w(n)$.) This part of Theorem 1 is used in the proof of part (v) of Theorem A but is not required for the proof of parts (i)-(iv). The notion of time-bounded Kolmogorov complexity can be informally explained as follows. Consider a (possibly algorithmically unsolvable) decision problem. Assume now that the Turing machine is allowed to use oracle information but it must solve the decision problem for any instance of size $\leq N$ in a time not exceeding a given computable function $t(N)$. Of course, the amount of oracle information used by the Turing machine can also depend on $N$. Here we assume that there is a natural notion of the size of instances of the decision problem. The oracle information is represented as a string of zeroes and ones. The amount of the oracle information is just the length of this sequence. Instead of the Turing machine one can imagine the computer program, say in FORTRAN or C, but using only the integer type of data. The program is supposed to work as follows. For every $N$ it recieves a string of 0's and 1's (the oracle information). After that it must be able to solve the decision problem for any input of size$\leq N$ in time bounded by $t(N)$. If $t$ is not too slowly growing, one can, for example, demand the list of all answers for all inputs of size $\leq N$ as the oracle information. (We assume that the number of instances of the problem of size $\leq N$ is finite for any $N$.) Roughly speaking, the time-bounded Kolmogorov complexity of the decision problem with time resources bounded by $t$ is the minimal amount of oracle information required to solve the decision problem for all inputs of size $\leq N$. The rigorous definition and discussion of the time-bounded Kolmogorov complexity can be found in [N2], [N4], [LV], [D]. According to [B] there exists a Turing machine $T$ such that its halting problem has time-bounded Kolmogorov complexity with time resources bounded by $t$ not less than $2^N/const(t)$ for any computable function $t$. Here inputs are finite strings of 0's and 1's, and the size of instances of the problem is just the length of the input strings (see also section 2 of [N2] and Theorem 2.5 in [ZL]). (In other words, there is no essentially better way to ask for the oracle information than just to ask for the list of all answers as in the example mentioned above. We cannot reduce the required amount of oracle information more than by a constant factor.) This theorem can be regarded as a quantitative version of Turing's classical theorem on the algorithmic unsolvability of the halting problem for Turing machines. In [N2] it was shown that the Barzdin theorem implies similar exponential lower bounds for the time-bounded Kolmogorov complexity of the word problem and the triviality problem for finitely presented groups as well as for the homeomorphism problem for compact $n$-dimensional manifolds, $n \geq 5$. These results are used in the proof of Theorem 1 and Theorem A(v) of the present paper.

The last part of Theorem 1 contains the same results for any fixed compact manifold of dimension $\geq 5$ (instead of $S^n$). (This almost immediately follows from the proof in the $S^n$ case.)

Very informally speaking, Theorem 1 implies that the space of all Riemannian structures on any compact $n$-dimensional manifold, $(n \geq 5)$, and, in particular, on $S^n$ has "many" "large" regions where it "looks" like the space of Riemannian structures on a manifold admitting a non-negatively curved Riemannian metric, and which can be reached starting from the standard metric only by a path passing through Riemannian metrics of extremely high curvature-pinching $(\sup |K|) \, diam^2$. (Most probably, it can be reached only by a path passing through Riemannian metrics of extremely high $-\inf K \, diam^2$, but at the moment we are not able to prove this.)

In section 1 we apply Theorem 1 to demonstrate that the graph of scalar curvature on the space of Riemannian metrics on $S^n$, $n \geq 5$, of *negative* constant scalar curvature has complicated geometry (see Theorem 2). In section 1 we also describe an approach to the construction of Einstein or almost Einstein metrics of *negative* scalar curvature on $S^n$ (Proposition 3). This approach can be used to reduce the problem of construction of Einstein or almost Einstein metrics of negative scalar curvature on $S^n$ to the problem of construction of such metrics on smooth homology spheres not admitting Riemannian metrics of non-negative scalar curvature, not representable as a non-trivial connected sum and such that the volume of these homology spheres endowed with any Riemannian metric of Ricci curvature $\geq -(n-1)$ is not less than a certain positive constant $v(n)$ depending only on $n$. But for manifolds with these properties one can hope to prove the existence of at least almost Einstein metrics by minimizing the volume over the space of Riemannian metrics of constant scalar curvature $-1$. (However, note that one will probably need an irreducibility condition much stronger than just impossibility to represent the homology sphere as a connected sum in a non-trivial way. Still, it seems quite possible that such extra conditions can be satisfied in the framework of our general approach.) Our technique can be similarly applied to an extremely wide class of systems of "equations" and "inequalities" involving Ricci, sectional and scalar curvature, volume, diameter, algebraic operations and $C^k-$ and $L^k-$norms (for any $k \geq 0$) (see Theorem C at the end of Section 1).

Theorem A follows from Theorem 1 by contradiction: Assuming Theorem A to be false, we construct (Proposition 0.2) an algorithm of the sort Theorem 1 precludes. This argument is given in section 4. For example, a very particular case of Proposition 0.2 states that if for all sufficiently large $x$ either $diam^{-1}((0,x]) \in R_1(S^n)$ is connected or there are no connected components on which the infimum of volume is $\geq v(n) > 0$ (or even if the number of connected components with this property grows with $x$ slower than $\exp(cx^n)$ for any $c > 0$), then there is an algorithm which distinguishes $S^n$ from any nonsimply-connected homology sphere with the property that its volume when endowed with any Riemannian metric of Ricci curvature $\geq -(n-1)$ is not less than $v(n)$. In order to explain here some ideas of our proof of Proposition 0.2 in section 4 we provide here the following argument illustrating possible connections between connectedness and algorithmic decidability (as well as between the number of connected components and the time-bounded Kolmogorov complexity).

Let $(X_k, d_k)$ be a family of precompact metric spaces, and $A_k \subset X_k$ be their subsets such that for any $a \in A_k$ and $x \in X_k \setminus A_k$ $d_k(a,x) \geq \theta(k) > 0$, where $\theta : \mathbb{N} \longrightarrow \mathbb{Q}$ is a computable function. Assume that for any $k$ $A_k$ is connected. Then one can try to solve the algorithmic problem "For given $k$ and $x \in X_k$ decide whether or not $x$ is an element of $A_k$" using the following algorithm: Construct a finite $\theta(k)/5$-net $N_k$ in $X_k$. (Of course, some extra assumptions about the family $\{(X_k, d_k)\}$ must be made in order to ensure the existence of an algorithm constructing for any $k$ such a net.) Construct a graph $Gr(k)$ such that the set of its vertices is in a bijective correspondence with the set of points of $N_k$. Two vertices are connected by an edge iff the distance between the corresponding points of $N_k$ is at most $2\theta(k)/3$. (Here we assume that the distances between points of $N_k$ can be approximately calculated. We perfom the computations of distances within to accuracy $\theta(k)/5$ and neglect the possible errors of computations.) Assume that we can find for any given $k$ $a_k \in A_k \bigcap N_k$. Now for any given $x \in X_k$ we can first find $y \in N_k$ $\theta(k)/2$-close to $x$ and then determine if the vertices of $Gr(k)$ corresponding to $y$ and $a_k$ are in the same component of $Gr(k)$. The answer will be positive if and only if $x \in A_k$. The connectedness of $A_k$ is, of course, crucial here.

Assume now that $A_k$ is not necessarily connected but we are allowed to use oracle information (depending on $k$) to solve the considered algorithmic problem. A reasonable way to do that is to require a description of a representative from every connected component of $A_k$ which must be also an element of $N_k$. Indeed, our assumption imply that for any component of $Gr(k)$ either all its vertices correspond to points in $A_k$ or none of them correspond to a point of $A_k$. Using the oracle information we can construct all components of $Gr(k)$ corresponding to points of $A_k$. Now for any given $x \in X_k$ find, as before, $y \in N_k$ $\theta(k)/2$-close to $x$ and check whether or not the corresponding vertex of $Gr(k)$ is in one of the constructed components. With some luck this algorithm will run in a time bounded by a computable function of $k$. The required amount of oracle information grows linearly with the number of connected components of $A_k$. Therefore any upper bound for the number of connected components of $A_k$ implies also an upper bound for the time-bounded Kolmogorov complexity of the considered decision problem.

Very roughly speaking, in section 4 we apply similar ideas in the situation when we have a continuous parameter $x$ instead of the discrete parameter $k$, $X_k$ (or, more precisely, $X_x$) is the set of all isometry classes of compact $n$-dimensional Riemannian manifolds such that $|K| \leq 1$, volume is $\geq v(n)/100$ and diameter is $\leq x$, and $A_k$ is the subset of $X_k$ formed by all Riemannian structures on a fixed manifold. The most difficult part of the proof is the construction of the required net on $X_x$. In fact, we construct this net in a neighborhood of $X_x$ using the following idea: Applying first the Ricci flow (cf. [Ha], [BMR]), then choosing appropriate harmonic local coordinates (cf. [JK]) and performing an algebraic approximation we can approximate any Riemannian manifold from $X_x$ by a Riemannian manifold from a compact finite-dimensional set. Then we construct a net in this finite-dimensional set.

Note that Theorems A, B are false as stated for $n = 2$ or 3. (This fact easily follows from the algorithmic recognizability of 2-dimensional manifolds and $S^3$ (cf. e.g. [T]).) However, we do not know whether or not a slightly weakened version of Theorems A,B where $\phi(x) = x$ is true. Moreover, we cannot say anything about the local minima of $diam$ on $R_1(M)$ in the case when $dim(M) = 2$ or 3 except the trivial remark that if $M$ is non-collapsible then the global minimum of $diam$ on $R_1(M)$ exists. On the other hand we believe that Theorems A and B are true for all compact 4-dimensional manifolds. It is not difficult to prove these theorems for manifolds $M$ such that $M = N^4 \# H^4 \# k(S^2 \times S^2)$, where $N^4$ is an arbitrary compact four-dimensional manifold, $H^4$ is an arbitrary compact four-dimensional manifold of non-zero simplicial volume (see [Gr 1] for the definition and properties of simplicial volume), and $k$ is the same absolute constant as in the text of Theorem 1.1B in [N2]. Indeed, the proof of Theorem A, (i)-(iv), and of Theorem B minus the exponential lower bound for the number of the deep local minima can be essentially found in [N3]. (These proofs can be obtained by combining the proofs of Theorems 9 and 11 in section 5 of [N3] with the remark at the beginning of section 5.A of [N3]). To prove Theorem A (v) and the exponential lower bound for the number of the deep local minima one can follow the proof of these statements given in the present paper but replacing Theorem 1 by the smooth version of Lemma 3.1(b) of [N2].

In a sequel of this paper we hope to prove the analogs of Theorem A and B in the situation when the diameter is regarded as a functional on the closure in $Met(M^n)$ of the space of Riemannian metrics on $M^n$ of sectional curvature $\geq -1$. (The elements of this closure are Alexandrov structures of curvature $\geq -1$ on $M^n$.)

## 1. Main algorithmic unsolvability result and its implications.

**Theorem 1. A.** For any fixed $n \geq 5$ there exists a positive rational $w(n)$ such that there is no algorithm which distinguishes between the standard $n$-dimensional sphere and smooth compact non-simply connected $n$-dimensional manifolds with the following properties: (1) The manifold does not admit a Riemannian metric of non-negative scalar curvature; (2) For any Riemannian metric on this manifold such that its Ricci curvature is greater than $-(n-1)$, the volume of the manifold with respect to this Riemannian metric is not less that $w(n)$; (3) The manifold cannot be represented as the connected sum of two smooth manifolds such that none of them is homeomorphic to the standard sphere.

**B.** More precisely, there exists an algorithm which constructs for any given $i \in \mathbb{N}$ a smooth $n$-dimensional homology sphere $S_i$ endowed with a Riemannian metric $\mu_i$ such that $S_i$ is either diffeomorphic to $S^n$ or is a non-simply connected homology sphere satisfying conditions (1)-(3) in the text of part A, and the following conditions are satisfied: (i) Denote by $I \subset \mathbb{N}$ the set of indices $i$ such thar $S_i$ is diffeomorphic to $S^n$. Then $I$ is non-recursive, and, moreover, for any Turing computable increasing function $\lambda$ the time-bounded Kolmogorov complexity $K^{(\lambda)}(I,N)$ of the decision problem "Is a given $i$ an element of $I$?" for all natural numbers $i \leq N$ with the time-bound $\lambda$ satisfies the following inequality:

$$K^{(\lambda)}(I,N) \geq N/c(\lambda,n) - const.$$

Here $c(\lambda,n)$ is a positive constant depending only on $\lambda$ and $n$ but not on $N$, and $const$ does not depend on $N$; (ii) For any $i$ the absolute values of sectional curvatures of $(S_i, \mu_i)$ do not exceed one, the convexity radius is not less than one, the diameter does not exceed $const_1(n)(\ln(i+1))^{\frac{1}{n}}$ and the volume is between $const_2(n)\ln(i+1)$ and $const_3(n)\ln(i+1)$, where $const_1(n), const_2(n)$ and $const_3(n)$ are some absolute positive constants.

**C.** Parts A and B of this theorem will be true for any compact manifold $M_0^n$ instead of $S^n$ providing that the following changes are made: First, one drops property (3). In the case if $M_0^n$ is not spin one drops also property (1). (In this last case only property (2) remains.) In part B $S_i$ must be smooth homology $M_0^n$ and not homology spheres.

**Remarks.** Part A of Theorem 1 is sufficient to prove Theorem A, parts (i)-(iv) and, as a corollary, a weaker version of Theorem B stated in the Remark after the text of Theorem B. We refer the reader to [N3], section 2 and [N4] for the definition and discussion of the time-bounded Kolmogorov complexity and the Barzdin theorem which is an essential ingredient in the proof of the part B of Theorem 1 (see also [LV], [D], [ZL]).

**Corollary:** For any $n \geq 5$ there exist $n$-dimensional smooth homology spheres satisfying conditions (1)-(3) in the text of the part A of Theorem 1.

**Remark.** The proof of the corollary is, in fact, a part of the proof of Theorem 1. From this point of view it would be more appropriate to state it as a theorem predcessing Theorem 1.

Formally speaking, the term "algorithm" is used in this paper as a synonym of the term "Turing machine". Also, to make Theorem 1 rigorous one must postulate how the (diffeomorphism types of) smooth manifolds and the (isometry classes) of Riemannian metrics are presented in a finite form (for computational purposes). One possibility to code manifolds in a finite form is Nash (i.e. smooth semialgebraic) atlases described in [BHP] (see also [N1] and the last section of the present paper). We can assume that all semialgebraic functions involved are defined over the field of real algebraic numbers. The Riemannian metrics $\mu_i$ in Theorem 1.B and its generalization in Theorem 1.C can be given by smooth semialgebraic (over the fields of real algebraic numbers) functions in the Nash local coordinates on the constructed manifolds. Another possibility is to use the Tognoli theorem stating that any smooth compact manifold is diffeomorphic to an algebraic subvariety of an Euclidean space. Using the Tarski-Seidenberg theorem it is not difficult to see that one can choose this algebraic subvariety as a zero set of a polynomial with algebraic coefficients (cf. [CtSh]; see also [BT]). One can regard this vector of coefficients as a finite set of data representing (the diffeomorphism type of) the considered manifold. Now the rigorous statement of Theorem 1.A is that for any fixed $n \geq 5$ there is no integer-valued recursive function defined on the set of codes of compact smooth $n$-dimensional manifolds (e.g. vectors of algebraic coefficients of $m$ polynomials of $n + m$ real variables such that their set of common zeroes is a compact $n$-dimensional manifold), assuming the value zero if the given manifold is diffeomorphic to $S^n$ and assuming the value one if the given manifold is a homology sphere which satsfies properties (1)-(3) stated in the text of Theorem 1.A. (We do not postulate what are the values of this function in the other cases.) This approach to coding manifolds can be similarly used in the situation of Theorem 1.B,C. In this approach we can choose Riemannian metrics induced by the embedding into the ambient Euclidean space as $\mu_i$ (of course, a certain care is necessary when the embeddings are constructed). In fact, we will prove that $S_i$ can be constructed as (non-singular algebraic) hypersurfaces in $\mathbb{R}^{n+1}$. So in the situation of Theorems 1.A, 1.B (but not Theorem 1.C!) we can assume that manifolds $S_i$ are coded by vectors of real algebraic coefficients of a polynomial on $\mathbb{R}^{n+1}$ such that $S_i$ is its zero set. (One can also show that it is possible to construct $S_i$ as non-singular algebraic hypersurfaces in $\mathbb{R}^{n+1}$ such that one can choose the Riemannian metric on $S_i$ induced by the embedding into $\mathbb{R}^{n+1}$ as $\mu_i$, but we do not need this fact here.)

We are going to show how to deduce Theorem 1 from Proposition 4 in Section 2 and to prove Proposition 4, thus completing the proof of Theorem 1 in Section 3. Proposition 4 is an algebraic reason d'être of smooth homology spheres satisfying the conditions (1)-(3) of the part A of Theorem 1. In the rest of this section we are going to describe some possible applications of Theorem 1 in Riemannian geometry different from Theorems A and B stated in the introduction. First, observe that a weaker and much simpler version of Theorem 1.A where the property (2) is not required implies that:

**Theorem 2. (Imprecise statement): A.** For any $n \geq 5$ and $k \geq 0$ there is no "nice" flow (e.g. no "nice" modification of the Ricci flow) on the space of Riemannian metrics on $S^n$ with the following property: The trajectory of the flow which starts at a Riemannian metric which is not of positive scalar curvature always reaches the set of Riemannian metrics of positive scalar curvature in time bounded by a computable function of volume, diameter and $C^k$-norm of the curvature tensor of the initial Riemannian metric. The adjective "nice" means here that the flow satisfies the following properties: 1) There exists an algorithm allowing one to trace (approximately in the $C^2$-norm) the trajectory of the flow starting from a sufficiently close approximation (in $C^4$-norm) to the initial Riemannian metric and until it reaches a Riemannian metric of positive scalar curvature; 2) The time of work of this algorithm can be majorized by a computable function of accuracy and metric invariants of the original metric (such as $C^k$-norm of the curvature tensor, injectivity radius, volume, diameter, etc.); 3) The algorithm can be applied to trace the trajectory of the flow on other compact $n$-dimensional manifolds (not necessarily $S^n$) to a given accuracy and during a given time unless by some reason it cannot continue its work. In this last case the algorithm informs that it cannot continue and stops. The initial data for this algorithm consist of a Nash atlas (in the sense of [BHP], see also [N1]) defined over the field of real algebraic numbers and a Nash (i.e. smooth semialgebraic for any local coordinate system) Riemannian metric (also defined over the field of real algebraic numbers).

**B.** There is no "nice" flow on the space of Riemannian metrics of constant scalar curvature and volume one on $S^n$, $n \geq 5$, with the following property: The trajectory of the flow which starts at a Riemannian metric of negative scalar curvature always reaches the set of Riemannian metrics of constant positive scalar curvature in a time bounded by a computable function of volume, diameter and the $C^k$-norm of the curvature tensor for some fixed $k$.

We do not want to obscure this principle by giving a more formal definition of "nice" flows. It is clear that this can be done.

The proof of A consists in demonstrating that the existence of such "nice" flow implies the existence of an algorithm distinguishing between the standard $n$-dimensional sphere and a smooth $n$-dimensional homology sphere not admitting a metric of non-negative scalar curvature. Indeed, consider the unknown manifold with the Riemannian metric induced by the embedding in the ambient Euclidean space. Compute an upper bound for the time of work of the algorithm tracing the flow and trace the flow during the computed time. If the flow will lead to a Riemannian metric of non-negative scalar curvature, then the manifold is $S^n$. If either the algorithm will stop prematurely or the flow will not reach a Riemannian metric of positive scalar curvature in the considered time, then the manifold is not $S^n$.

The proof of part B is obtained by a quite similar argument after one takes into consideration the fact that there exists exactly one Riemannian metric of constant negative scalar curvature in the conformal class of any Riemannian metric of (not necessarily constant) negative scalar curvature (cf. [Sch] and references there), and that this constant scalar curvature metric is the solution of a semilinear elliptic equation which can be numerically solved with any prescribed accuracy.

*Thus, it is completely wrong to visualize the graph of the scalar curvature on the space of Riemannian structures on $S^n$ of constant non-positive scalar curvature and volume*

*one as a "mountain slope" evenly raising (without any "peaks" or "precipices") from the domain of large negative values of scalar curvature towards the domain on the graph, where the scalar curvature is positive.*

Let us describe now another possible application of Theorem 1. In [N1] we observed that the non-existence of an algorithm recognizing $S^n$, ($n \geq 5$), implies that if any homology sphere with infinite fundamental group almost admits an Einstein metric then the standard sphere $S^n$ almost admits an Einstein metric of negative scalar curvature. We say that a compact manifold $M$ *almost admits an Einstein metric of negative scalar curvature* if there exists a sequence of Riemannian metrics $\{g_i\}$ on $M$ of constant negative scalar curvature $s_i$ and of diameter one such that: (0) The sequence $s_i$ is uniformly bounded from below; (i)
$$\lim_{i \longrightarrow \infty} \|Ric(g_i) - (s_i/n)g_i\|_{C^0} = 0,$$
where the norms are taken with respect to $g_i$; (ii) The sequence of volumes of $(M, g_i)$ is bounded from below by a positive constant. As in [N1], if it is desirable, then one can modify this definition in many ways. For example, one can replace the $C^0$ norm in (i) by any $C^k$ or $L^k$ norm and/or choose any non-negative integer $k$ and any computable increasing function of the $C^k$-norm of the curvature tensor $R$ (say, $\exp(\|R\|_{C^2})$) and to demand instead of (i) that the product of the value of this function and $\|Ric(g_i) - (s_i/n)g_i\|_{C^0}$ must tend to zero. (All results below will remain valid, if one changes the definition in such a way. "Computable" means here that the restriction of this function on integers is (Turing) computable.) Moreover, it was proven in [N1] that if the very plausible conjecture that there exists an algorithm solving a problem which can be informally stated as the question "Is a given almost Einstein metric close to an Einstein metric?", then the existence of Einstein metric on all smooth homology spheres of a fixed dimension $n \geq 5$ with infinite fundamental group implies the existence of infinitely many different Einstein structures of volume one and *negative* scalar curvature on the sphere $S^n$. (See [N1] for the rigorous statement of the mentioned algorithmic problem. It is extremely plausible that this algorithmic problem is algorithmically solvable.) Thus, one can try to prove the existence of (almost) Einstein metrics of negative scalar curvature on $S^n$ by proving the existence of (almost) Einstein metrics on all smooth $n$-dimensional homology spheres with infinite fundamental group. However, there is no known way to do that. One general method which could be used to show that a homology sphere $H$ admits an Einstein metric is the following: Consider the space of Riemannian structures of volume one and constant scalar curvature on $H$. The scalar curvature can be regarded as the functional on this space. One can try to maximize this functional. It is known the the supremum of this functional does not exceed the scalar curvature of the standard metric of volume one $S^n$. Further, if there exists a Riemannian structure on $H$ at which this maximum is achieved and the corresponding maximal value is non-positive then this Riemannian structure is Einstein. (The proofs of these results can be found in [Sch] (see also [An2]).) The last result indicates that ad hoc it might be simpler to prove that the described procedure leads to an Einstein (or almost Einstein) metric if the manifold does not admit Riemannian metrics of non-negative scalar curvature. (There are some additional technical reasons why it might be simpler to consider only manifolds not admitting Riemannian metrics of non-negative scalar curvature. For example, for such manifolds there exists a unique metric of constant scalar curvature in every conformal class, and the space of Riemannian metrics of constant scalar curvature and volume one is known to be a ILH-manifold (see [Be]).) As it is clear from the work of Anderson ([An 2]) about Einstein metrics on compact 3-dimensional manifolds, the impossibility to represent a manifold as the connected sum of two manifolds none of which is homeomorphic to the sphere also can help to prove the existence of an (almost) Einstein metric on this manifold. Theorem 1 allows us to replace the problem of finding (almost) Einstein metrics on all homology spheres by the problem of finding such metrics only on homology spheres not admitting metrics of non-negative scalar curvature and such that they cannot be represented as the connected sum in a non-trivial way. An equivalent dual version of the just desrcibed approach to construction of Einstein metrics is to minimize volume over the space of metrics of constant negative scalar curvature $-n$. This seems to be easier if the volume of the manifold with respect to any Riemannian metric of Ricci curvature $\geq -(n-1)$ is greater than $v(n) > 0$ (as it happens, for example, for manifolds of non-zero simplicial volume; see [Gr 1]). In this case the condition (ii) in the definition of almost Einstein metrics of negative scalar curvature automatically follows from conditions (0) and (i). All these remarks indicate that the following proposition can be an important step in the proof that $S^n$, $n \geq 5$, almost admits Einstein metrics of negative scalar curvature.

**Proposition 3.** Let $n \geq 5$ be fixed. Assume that any smooth $n$-dimensional homology sphere satsfying properties (1)-(3) from the text of Theorem 1.A almost admits an Einstein metric. Then the standard sphere $S^n$ almost admits an Einstein metric of negative scalar curvature.

**Remark.** As in [N1] one can also state the version of this proposition for Einstein metrics instead of almost Einstein metrics providing that either the conjecture (D) on p. 78 of [N1] is true, or that one assumes not merely the existence of Einstein metrics on homology spheres satisfying conditions (1)-(3) from the text of Theorem 1.A but the existence of *stable* Einstein metrics in the sense of [Be] (Definition 4.63). As in [N1] one can prove under this assumption the existence of infinitely many distinct Einstein structures of scalar curvature $-1$ on $S^n$, ($n \geq 5$).

**Sketch of the proof.** This theorem follows from Theorem 1 exactly as Theorem 1 of [N1] followed from the S. Novikov theorem stating the algorithmic unrecognizability of $S^n$, $n \geq 5$. The idea of the proof by contradiction is to assume that any smooth $n$-dimensional homology sphere satisfying properties (1)-(3) from the text of Theorem 1.A almost admits an Einstein metric and that $S^n$ does not almost admit an Einstein metric of negative scalar curvature and to construct an algorithm which distinguishes between $S^n$ and an irreducible homology sphere. The idea of the algorithm is very simple: We look (using a trial and error algorithm) for either a metric of positive Ricci curvature on a given manifold or an almost Einstein metric. It follows from our assumptions that we eventually will find either a metric of positive Ricci curvature, and in this case the manifold is $S^n$, or an almost Einstein metric, and we know that in this last case the manifold cannot be $S^n$. □

Of course, the same idea to combine Theorem 1 with [N1] can be applied not only to find (almost) Einstein metrics but when one tries to prove the existence of an arbitrary finite systems of equations and inequalities on the space of Riemannian metrics on a fixed compact manifold of dimension $\geq 5$. Assume that some such system involving Riemannian metric, curvature tensor and/or its separate components, sectional curvature, Ricci curvature, scalar curvature, injectivity radius, convexity radius, volume, diameter, $C^k$, $L^k$ and $W^{k,p}$-norms, algebraic operations and, more generally, elementary functions be given. (Example: $Ric - (scal/n)g = 0$ and $scal < 0$.) Assume that any equation or inequality is written in the form $E?0$, where $E$ is a Riemannian functional and ? is $=, >, <, \geq$ or $\leq$. For any such system $S$, any fixed $k$, and any $\epsilon$ we can replace $S$ by a one-parametric family of systems $S_\epsilon$ indexed by positive $\epsilon$ and obtained from $S$ as follows: Any equation $E = 0$ is replaced by $\|E\|_{C^k} < \epsilon$, all strict inequalities are preserved as they are, any inequality $E \leq 0$ (or $E \geq 0$) is replaced by $E < \epsilon$ (or, correspondingly, $E > -\epsilon$).

**Theorem C.** For any $n \geq 5$ and any compact manifold $M^n$ of dimension $n$ consider the class $\sigma$ of all smooth $n$-dimensional manifolds with the same homology as $M^n$ and such that the volume of any of these manifolds endowed with any Riemannian metric on this manifold of Ricci curvature greater than $-(n-1)$ is not less than the positive constant $w(n)$ depending only on $n$ defined in the text of Theorem 1.A. If $M^n$ is spin we impose an additional constraint in the definition of $\sigma$ leaving only those manifolds which do not admit a Riemannian metric of non-negative scalar curvature. If for any $\epsilon > 0$ the system $S_\epsilon$ has a solution in the space of Riemannian metrics on any manifold from the class $\sigma$ then for any $\epsilon > 0$ $S_\epsilon$ has a solution in the space of Riemannian metrics on $M^n$.

## 2. Reduction of Theorem 1 to homological group theory.

In this section we are going to deduce Theorem 1 from the following proposition which will be proven in the next section:

**Proposition 4.** For any $n \geq 5$ there exists a group $A_0$ and a homology class $h_n \in H_n(BA_0)$ such that:
(1) $A_0$ belongs to the class $P$ of finitely presented groups defined as follows: $P$ is the minimal class which contains all discrete subgroups of connected Lie groups and is closed with respect to forming amalgamated free products of groups from this class such that the amalgamated subgroup is free or abelian and forming HNN-extensions of groups from this class where the associated subgroup is free or abelian;
(2) $H_1(A_0) = \{0\}$;
(3) Let $h_n^+$ denote the image of $h_n$ in $H_n(BA_0^+)$ under the isomorphism of the $n$th homology groups of $BA_0$ and $BA_0^+$, where $BA_0^+$ denotes the result of the application of the Quillen +-construction to $BA_0$. Then the class $h_n^+$ is spherical (i.e. it is in the image of the Hurewicz homomorphism $\pi_n(BA_0^+) \longrightarrow H_n(BA_0^+)$.);
(4) For any non-zero integer $k$ $kh_n \neq 0$ in $H_n(BA_0)$;
and either
(5) The simplicial norm of $h_n$ regarded as a real homology class is non-zero;
or
(5'.1) The class $h_n$ is not in the image of the homomorphism $p_* : H_n(Q_{A_0}//A_0) \longrightarrow H_n(A_0)$ defined on the p. 79 of [Gr 2]; and
(5'.2) The group $A_0$ is torsion-free, and there exist elements $g_1, \ldots, g_l$ normally generating

$A_0$ and satisfying the following additional condition: None of them is a square of an element of $A_0$.

The Quillen +-contruction and its properties are explained, for example, in [Ros 2]. In particular, it assigns to any CW-compex $X$ and a perfect normal subgroup $\pi$ of $\pi_1(X)$ a CW-complex $X^+$ and the canonical inclusion $X \longrightarrow X^+$ inducing the epimorphism of fundamental groups $\pi_1(X) \longrightarrow \pi_1(X^+) = \pi_1(X)/\pi$ and such that the pair $(X^+, X)$ is homologically acyclic. $X^+$ can be obtained from $X$ by first attaching 2-cells in order to kill $\pi$ and then attaching 3-cells to restore the second homology groups of $X$. In the present paper we apply the Quillen +-construction only in the case when $[\pi_1(X), \pi_1(X)]$ is perfect and $\pi = [\pi_1(X), \pi_1(X)]$. Therefore from now on we will be omitting the mention of $\pi$ while discussing the +-construction. Moreover, in the most cases we will be applying the +-construction in the situation when $\pi_1(X)$ is perfect and, thus, $\pi = [\pi_1(X), \pi_1(X)] = \pi_1(X)$. In this case $X^+$ is simply-connected and the canonical inclusion induces the isomorphism of homology groups of $X$ and $X^+$. Finally, recall that the +-construction is functorial. The Quillen +-construction appears in our work because of the following result of J. Hausmann and P. Vogel (cf. [Haus]) based on the theory of homology surgery by S. Cappell and J. Shaneson ([CS]): Let $n \geq 5$ be any fixed number, $A$ be a finitely presented group and $h$ be an element of $H_n(A)$. The necessary and sufficient conditions of the existence of a smooth $n$-dimensional homology sphere $\Sigma^n$ with the fundamental group $A$ such that $h$ is the image of the fundamental homology class of $\Sigma^n$ under the homomorphism $H_n(\Sigma^n) \longrightarrow H_n(BA) = H_n(A)$ induced by the classifying map into $BA$ are the following:
(1) $H_1(A) = H_2(A) = \{0\}$;
(2) The image $h^+$ of $h$ in $H_n(BA^+)$ under the isomorphism induced by the canonical embedding $BA \longrightarrow BA^+$ is in the image of the Hurewicz homomorphism $\pi_n(BA^+) \longrightarrow H_n(BA^+)$.
It is a well-known result of Kervaire ([K0]) that (1) is the necessary and sufficient condition for a finitely-presented group $A$ to be the fundamental group of a smooth $n$-dimensional homology sphere (for any fixed $n \geq 5$). The necessity of (2) is obvious: If there exists a homology sphere $\Sigma^n$ such that $\pi_1(\Sigma^n) = A$ and $h$ is the image of the fundamental homology class of $\Sigma^n$ under the homomorphism of homology groups induced by the classifying map $f : \Sigma^n \longrightarrow BA$, then applying the Quillen +-construction to the map $f$ we obtain a continuous map $f^+ : (\Sigma^n)^+ = S^n \longrightarrow BA^+$. The induced homomorphism of homology groups maps the fundamental homology class of $S^n$ into $h^+$.

Following [Gr 1] one can define for any topological space $X$ the simplicial norm of any real homology class $z \in H_k(X, \mathbb{R})$ as the infimum over all singular chains $\Sigma_i \alpha_i \sigma_i$, ($\alpha_i \in \mathbb{R}$, $\sigma_i$ is a singular simplex in $X$), representing $z$ of the $l^1-$norm of $z$ $\|z\| = \Sigma_i |\alpha_i|$. For a compact orientable manifold $M$ the simplicial volume $\|M\|$ of $M$ is defined as the simplicial norm of the real fundamental homology class of $M$. If $M$ is a compact non-orientable manifold, then its simplicial volume $\|M\|$ is defined as $\frac{1}{2}\|\tilde{M}\|$, where $\tilde{M}$ is the orientable double covering of $M$. The importance of this notion for our purposes is due to the fact that for any Riemannian manifold $(M, g)$ such that $Ric_{(M,g)} \geq -(n-1)$ the volume of $(M, g)$ is not less than $const(n)\|M\|$, where $const(n)$ is a positive constant which can be explicitly written down ([Gr 1], Main inequality on p. 12). Moreover, the simplicial

volume of an orientable manifold $M$ coincides with the simplicial norm of the image of the real fundamental homology class of $M$ in $H_n(K(\pi_1(M),1), \mathbb{R})$, ($n = dim\ M$), under the homomorphism induced by the classifying map $M \longrightarrow K(\pi_1(M),1)$ (Corollary (B) on p. 40 of [Gr 1]). Many other properties of simplicial volume and a closely related notion of bounded cohomology can be found in the seminal paper [Gr 1] (see also [I], [BP]).

Furthermore, $Q_{A_0}$ denotes the simplicial complex of almost nilpotent groups of $A_0$. (Its vertices are in a bijective correspondance with all nontrivial elements of $A_0$; the simplex formed by several vertices is a simplex of $Q_{A_0}$ if and only if the corresponding elements of $A_0$ generate an almost nilpotent subgroup of $A_0$.) There is a natural action of $A_0$ on $Q_{A_0}$ induced by the action of $A_0$ on itself by conjugation. Thus, we can consider the Borel construction $p : (Q_{A_0})_{A_0} \equiv Q_{A_0}//A_0 \longrightarrow BA_0$. (See either chapter 6.6 of [Gr 2] or the beginning of the proof of Lemma 5.8 below for a more detailed description of $Q_{A_0}//A_0$ and $p$.) Now we can state another result of Gromov used in our paper: Assume that $G$ is a discrete group such that the homomorphism $p_* : H_n(Q_G//G) \longrightarrow H_n(G)$ induced by $p$ vanishes for some $n$. Assume that $M$ is a compact $n$-dimensional manifold, and $q : M \longrightarrow BG$ is a continuous map such that the image of the fundamental homology class of $M$ under the homomorphism $q_* : H_n(M) \longrightarrow H_n(BG)$ induced by $q$ is not equal to zero. Then Theorem 6.6.D of [Gr 2] implies (as in the proof of Theorem 6.6.D" in [Gr 2]) that there exists a positive $c(n)$ such that $vol(M,g) \geq c(n)$ for any Riemannian metric $g$ on $M$ such that $Ric_{(M,g)} \geq -(n-1)$. (In fact, it is clear from the proof of Theorem 6.6.D in [Gr 2] that instead of vanishing of $p_*$ on $H_n(Q_G//G)$ one can demand a weaker condition that the image of the fundamental homology class of $M$ under the homomorphism $q_*$ is not in the range of $p_*$.)

Below we give two parallel but different proofs of Theorem 1 (and, thus, of Theorems A,B,C). One proof is based on [Gr 1] and on constructions in section 3 below of $A_0$ satisfying conditions (1)-(5). Another proof uses [Gr 2] and constructions of groups $A_0$ satysfying (1)-(4), (5'.1),(5'.2) also described in section 3.

Here is the scheme of the proof of Theorem 1 assuming Proposition 4. Start from a Turing machine $T$ with the unsolvable halting problem. First we are going to effectively construct a sequence of (finite presentations of) groups $G_i$ and homology classes $\bar{h}_i \in H_n(G_i)$ such that
(i) The sequence of indices $i$ for which $G_i$ is non-trivial is recursively enumerable but non-recursive (that is, there is no algorithm deciding for a given $i$ whether or not $G_i$ is trivial). Moreover, the halting problem for $T$ for inputs of length $\leq N$ is effectively reducible to the triviality problem for the first $2^N$ groups $G_i$. The numbers of generators and relators in constructed finite presentations $G_i$ do not depend on $i$. The lengths of relators are bounded by a linear function of $\ln(i+1)$.
(ii) If $G_i$ is non-trivial, then for any integer non-zero $k$ $k\bar{h}_i \neq 0$;
(iii) $G_i$ is the universal central extension of a perfect group $L_i$, which belongs to the class $P$ introduced in the text of Proposition 4. (This implies, in particular, that $H_1(G_i) = H_2(G_i) = 0$);
(iv) The image $\bar{h}_i^+$ of $\bar{h}_i$ in $H_n(BG_i^+)$ is spherical;
(v) If $G_i$ is non-trivial, then either (v.i) the simplicial norm of $\bar{h}_i$ is a positive number $r$ not

depending on $i$ or (v.ii) $\bar{h}_i$ is not in the image of the homomorphism $p_* : H_n(Q_{L_i}//L_i) \longrightarrow H_n(L_i)$ defined as on p. 79 of [Gr 2];
and
(vi) $G_i$ is not a free product of non-trivial groups.

Then we effectively construct a sequence of smooth $n$-dimensional homology spheres $S_i$ such that (a) $\pi_1(S_i) = G_i$ and (b) If $G_i$ is non-trivial then the image of the fundamental homology class of $G_i$ under the classifying map in $H_n(BG_i)$ is $\bar{h}_i$. Note that the mentioned result of J. Hausmann and P. Vogel ([Haus]) and conditions (iii) and (iv) immediately imply the existence of smooth $n$-dimensional homology spheres with properties (a), (b).

Main Inequality on p. 12 of [Gr 1], Theorem 6.6 D and the proof of Theorem 6.6.D" in [Gr 2], and the property (v) of $G_i$ imply that if $S_i$ has non-trivial fundamental group (that is, $S_i$ is not homeomorphic to $S^n$), then for any Riemannian metric of Ricci curvatue $\geq -(n-1)$ on $S_i$ the volume of $S_i$ with respect to this metric is not less than some positive constant depending only on $n$ and $r$ if (v.i) holds and $n$ if (v.ii) holds. (Remark: More precisely, the assumptions of Theorem 6.6 D of [Gr 2] are stronger then (v.ii). Namely, in our situation the assumption of Theorem 6.6 D is that the homomorphism $p_*$ is trivial on $H_n(Q_{L_i}//L_i)$. But, as it was already noted, the proof of Theorem 6.6 D in [Gr 2] is valid under the weaker assumption that the image of the fundamental homology class of $V$ in $H_n(K(\pi,1))$ is not in the image of the homomophism $p_* : H_n(Q//\pi) \longrightarrow H_n(K(\pi,1))$. (Here we are using the notations of [Gr 2].) Although our proofs of (5.1') for specific groups $A_0$ in section 3 below imply also the triviality of the homomorphism $p_*$, we prefer this slightly stronger version of Theorem 6.6 D.) The impossibility to represent the homology spheres $S_i$ as connected sum of manifolds none of which is a homotopy sphere immediately follows from (vi). The condition (i) and classical result of the Gromov-Lawson-Rosenberg ([GL 2], [Ros 1]) imply that if $S_i$ is non-trivial then it does not admit a Riemannian metric of non-negative scalar curvature since groups from the class $P$ and their central extensions are known to satisfy the Baum-Connes conjecture that the Baum-Connes assembly map is a split injection, and, hence, the strong Novikov conjecture (SNC 2) in terminology of [Ros 1]. (A proof of this fact which is essentially a compilation of results of Baum-Connes ([BC], [BCH]), Pimsner ([Pi 1], [Pi 2]) and Kasparov-Skandalis ([KaS]) can be found in [BW]. A reader interested only in the proof of Theorems A and B and not willing to study the proof of Baum-Connes conjecture for groups from the class $P$ can just delete property (1) from the text of Theorem 1 above. The proof of Theorem A (and therefore of Theorem B) given here uses only this weaker version of Theorem 1.)

Thus, we obtain the algorithmic unsolvability of the problem of recognition of $S^n$ up to homeomorphism in the class of smooth homology spheres which are either homeomorphic to $S^n$ or do not admit Riemannian metrics of non-positive scalar curvature, do not admit Riemannian metrics of arbitrarily small volume and $Ric > -(n-1)$, and cannot be represented as the connected sum of two manifolds none of which is a homotopy sphere. Then we prove that the algorithmic unsolvability of the just stated homeomorphism problem implies the algorithmic unsolvability of the corresponding diffeomorphism problem, that is, Theorem 1.A. Alternatively, by a little change in the construction of $S_i$ we can make $S_i$ embeddable into $\mathbb{R}^{n+1}$ for any $i$. Since any smooth $n$-dimensional homotopy sphere embeddable into $\mathbb{R}^{n+1}$ is diffeomorphic to $S^n$, $S_i$ is homeomorphic to $S^n$ if and only if $S_i$ is diffeomorphic to $S^n$, and we again obtain Theorem 1.A. The next step is to introduce Riemannian metrics on the constructed homology spheres $S_i$ such that the convexity radius is at least one, the absolute values of sectional curvatures are between $-1$ and $1$, the diameter does not exceed $const_1(n)(\ln(i+1))^{\frac{1}{n}}$ and the volume is between $const_2(n)\ln(i+1)$ and $const_3(n)\ln(i+1)$. This part of the proof heavily relies on the material of section 2 of [N2]. Observe, that in the process of proving of Theorem 1.A we not only prove the algorithmic unsolvability of the considered diffeomorphism problem but obtain the reduction of the halting problem for an arbitrary fixed Turing machine to the considered diffeomorphism problem. Therefore we can use the Barzdin theorem stating the existence of a Turing machine $T_0$ such that the time-bounded Kolmogorov complexity of the halting problem for $T_0$ grows exponentially with the length of the inputs for any fixed computable bound on time. Starting our construction from the Turing machine $T = T_0$ we ensure that for any computable time bound the problem of recognition of $S^n$ in the constructed sequence of smooth Riemannian homology spheres has the time-bounded Kolmogorov complexity greater than an exponential function of the volume.

Finally, to generalize this argument for an arbitrary compact manifold $M_0^n$ (instead of $S^n$) we just take Riemannian connected sums of the constructed smooth Riemannian homology spheres with a copy of $M_0^n$ endowed with an arbitrary Riemannian metric. The only statement which is not immediately obvious is that the connected sum of $M_0^n$ and any of the constructed nonsimply-connected homology spheres not admitting a Riemannian metric of non-negative scalar curvature also does not admit a Riemannian metric of non-negative scalar curvature. In fact, we need the assumption that $M_0^n$ is spin in order to demonstrate this fact.

**Step 1. From $A_0$ to $G_i$**

To construct $G_i$ first observe that for any Turing machine $T$ and its input $\omega$ one can construct a finite presentation of a torsion-free group $G(T)$ depending only on $T$ and a word $w(T, \omega) \in G(T)$ such that $T$ halts with $\omega$ if and only if $w(T, \omega) = e$ in $G(T)$. (cf. [Rot]). (This fact implies the algorithmic unsolvability of the triviality problem for finitely presented groups.) The construction in [Rot] yields a group in the class $P$ defined in the text of Proposition 4. Let $T$ be any Turing machine such that its halting set is not recursive (equivalently, the halting problem for $T$ is unsolvable. We can identify the set of all inputs of $T$ with the set of positive integer numbers.) At this stage we do not want to impose any further restrictions on $T$. But on Step 5 we will make a particular choice of $T$ in order to prove the statement about the time-bounded Kolmogorov complexity in the part B of Theorem 1. Applying the "witness" construction as in [M], pp. 13-14 to $G(T)$ and words $w_i = w(T, i)$ for all $i$ we obtain a recursive sequence of finite presentations of perfect groups $\{W_i\}$ such that the set of indices $i$ such that $W_i$ is trivial is a non-recursive recursively enumerable set. $W_i$ has the following finite presentation: The set of generators of $W_i$ includes all generators $x_1, \ldots x_k$ of $G(T)$ plus three new generators $a, b, c$. Its set of relations includes all relations of $G(T)$ as well as the following new relations:

$$(1) \quad a^{-1}ba = c^{-1}b^{-1}cbc$$

$$(2) \quad a^{-2}b^{-1}aba^2 = c^{-2}b^{-1}cbc^2$$

$$(3) \quad a^{-3}[w_i,b]a^3 = c^{-3}bc^3$$

$$(4) \quad a^{-(3+j)}x_j b a^{(3+j)} = c^{-(3+j)}bc^{(3+j)}; j = 1,\ldots,k.$$

It is clear from this finite presentation that $H_1(W_i) = \{0\}$. It is not difficult to see ([M], p.13) that if $w_i \neq e$ in $G(T)$ $W_i$ is the amalgamated free product

$$(G(T)* < a,b_1| >)* < b_2,c| >; \quad (*)$$

$$V_1 = V_2, \quad (**)$$

where $V_1$ is the free groups generated by $b_1$ and the left hand sides of relations (1)-(4), where $b$ is replaced by $b_1$; $V_2$ is the free group generated by $b_2$ and the right hand sides of the relations (1)-(4), where $b$ is replaced by $b_2$. (We replaced $b$ by two different letters $b_1$, $b_2$ but added the equation identifying these elements.) It is easy to see that $W_i$ is the normal closure of the element $c$. It is clear that $c$ is an element of infinite order. The construction of groups $G(T)$ in [Rot] implies that the cohomological dimension of $W_i$ is at most two (cf. [Mil], Theorem 4.12), and that $W_i$ belong to the class $P$ defined in the text of Proposition 4. ($G(T)$ is the result of a sequence of three HNN-extensions with associated free groups.) Now let $g_1,\ldots,g_l$ be a set of normal generators of $A_0$. If $A_0$ satisfies conditions (5'.1), (5'.2), then we assume that these generators are the same as in (5'.2). In particular, all of them are of infinite order. In this case consider $l$ copies $W_i^{(1)},\ldots,W_i^{(l)}$ of $W_i$. Each of these copies $W_i^{(j)}$ is normally generated by a copy $c^{(j)}$ of the element $c$. Let $L_i = < A_0 * W_i^{(1)} * \ldots * W_i^{(l)} | g_1 = c^{(1)},\ldots,g_l = c^{(l)} >$. If $A_0$ is known to satisfy condition (5) but not necessarily (5'.2) then it is possible that some of the normal generators, say $g_{k+1},\ldots,g_l$ are of finite order. In this case consider for any $j \in \{k+1,\ldots,l\}$ the following group $V_i^{(j)}$. Take the free product of the cyclic subgroup $\mathbb{Z}_{ord(g_j)}$ and $G(T)$. $V_i^{(j)}$ is the result of the application of the "witness" construction described above to this group and its element $w_i = w(T,i)$ coming from $G(T)$. Exactly as before this group will be trivial if and only if $T$ halts with input $i$. The free product of the finite cyclic group and $G(T)$ embeds in $V_i^{(j)}$ if $V_i^{(j)}$ is non-trivial. Denote by $z_j$ the image of the generator of the cyclic group in $V_i^{(j)}$. The element $z_j$ has the same order as $g_j$. Let $L_i = < A_0 * W_i^{(1)} \ldots * W_i^{(k)} * V_i^{(k+1)} * \ldots * V_i^{(l)} | g_1 = c^{(1)},\ldots,g_k = c^{(k)}, g_{k+1} = z_{k+1},\ldots,g_l = z_l >$. It is clear that $L_i$ will be trivial if and only if $w_i \in G(T)$ is trivial. Hence the set of indices $i$ for which $L_i$ is trivial is non-recursive. Using the Mayer-Vietoris sequence one can easily see that $L_i$ is perfect. Denote by $q_i$ the image of $h_n \in H_n(A_0)$ in $H_n(L_i)$ under the inclusion homomorphism. Obviously $q_i = 0$ if and only if $w_i \in G(T)$ is trivial.

Let $G_i$ be the universal central extension of $L_i$. Many properties of the universal central extensions can be found in [K], [Ros 2], [Mil 1]. In particular, $G_i$ is trivial if and only if $L_i$ is trivial; $H_1(G_i) = 0$, and $H_2(G_i) = 0$. The universal central extension of a perfect group $L$ with the generators $f_i$, $i \in \{1,\ldots,p\}$, and relations $r_j = e$, $j \in \{1,\ldots,q\}$, can be described by the following finite presentation: The set of generators coincides with the set of generators of $L$. The set of relators consists of all commutators $[f_i, r_j]$, $i \in \{1, \ldots, p\}$, $q \in \{1, \ldots q\}$, and words $\lambda_i, i \in \{1, \ldots p\}$, satisfying the following condition: Denote by $F$ the free group generated by all generators of $L$ and by $R$ its subgroup generated by all relators $r_j$ of $L$. Then $\lambda_i \in R$ in $F$ (when regarded as an element of $F$) and $f_i = \lambda_i c_i$ in $F$, where $c_i \in [F, F]$. The existence of $\lambda_i$ follows from the perfectness of $L$. Also, if $C$ is any perfect group and $\bar{C}$ is its universal central extension then there is the following exact sequence:

$$0 \longrightarrow H_2(C) \longrightarrow \bar{C} \longrightarrow C \longrightarrow 1. \quad (***)$$

We claim that $q_i \in H_n(L_i)$ will be the image of some element $\bar{h}_i$ from $H_n(G_i)$ under the canonical homomorphism $H_n(G_i) \longrightarrow H_n(L_i)$ and, moreover this element can be chosen such that its image $\bar{h}_i^+$ in $H_n(BG_i^+)$ is spherical. Indeed, first note that for considered values of $n$ $H_n(BA_0)$ is isomorphic to $H_n(BL_i)$ if $w_i$ is non-trivial. In this case $H_n(BA_0^+)$ is isomorphic to $H_n(BL_i^+)$, and looking at the commutative square formed by $\pi_n(BA_0^+)$, $H_n(BA_0^+)$, $\pi_n(BL_i^+)$, $H_n(BL_i^+)$ with arrows induced by the inclusion $A_0 \longrightarrow L_i$ and the Hurewicz homomorphisms we see that the property (3) of the group $A_0$ postulated in Proposition 4 implies that the image $q_i^+$ of $q_i$ in $H_n(BL_i^+)$ is the image of some element $\hat{q}_i \in \pi_n(BL_i^+)$. Now note that for any perfect group $C$ and its universal central extension $\bar{C}$ the map $B\bar{C}^+ \longrightarrow BC^+$ is a fibration with the fiber $BH_2(C)$ ([Ros 2, Remark 5.2.9]). Therefore for the considered values of $n$ the homomorphism $\bar{C} \longrightarrow C$ induces the isomorphism of $\pi_n(B\bar{C}^+)$ and $\pi_n(BC^+)$. Let us apply this fact for $C = L_i$ and $\bar{C} = G_i$. Denote the inverse image of $\hat{q}_i$ in $\pi_n(BG_i^+)$ by $\check{q}_i$. Now we can define $\bar{h}_i^+$ as the image of $\check{q}_i$ under the Hurewicz homomorphism $\pi_n(BG_i^+) \longrightarrow H_n(BG_i^+)$. Define $\bar{h}_i$ as the inverse image of $\bar{h}_i^+$ under the isomorphism $H_n(BG_i) \longrightarrow H_n(BG_i^+)$. Thus, we constructed a sequence of groups $G_i$ and elements $\bar{h}_i \in H_n(BG_i)$ which has properties (i)-(iv). In order to see that if $A_0$ has property (5) and $G_i$ is non-trivial, then $G_i$ has property (v.i) we can first invoke Theorem on p. 55 of [Gr 1] which implies that in this case the class $q_i$ regarded as a real class has the same (non-zero) simplicial norm as $h_n$. Since $q_i$ is the image of $\bar{h}_i$ $\bar{h}_i$ also has a non-zero simplicial norm. Passing to the dual real cohomology classes, using the Mapping Theorem on p. 40 of [Gr 1] and returning to homology classes we can conclude that $\bar{h}_i$ has the same simplicial norm as $q_i$ and $h_n$.

Now assume that $A_0$ has properties (5'.1) and (5'.2), and $G_i$ is non-trivial. The postulated property (v.ii) of $G_i$ will be proven using the sequence of lemmae below. Let $A$ be any torsion-free finitely presented group, and $g \in A$ be any element such that for any $a \in A$ $g \neq a^2$. Let $B$ be obtained as the result of the application of the "witness" construction as in [M], pp. 13-14 to any finitely presented torsion-free group $X$ of cohomology dimension $\leq 2$ and its *non-trivial* element $w$. Denote generators of $X$ by $g_1, \ldots g_r$. $B$ has the following finite presentation: Its set of generators includes all generators $g_1, \ldots g_r$ of $X$ plus three new generators $a, b, c$. Its set of relations includes all relations of $X$ as well as the following new relations:

$$(5) \quad a^{-1}ba = c^{-1}b^{-1}cbc$$

$$(6) \quad a^{-2}b^{-1}aba^2 = c^{-2}b^{-1}cbc^2$$

$$(7) \quad a^{-3}[w,b]a^3 = c^{-3}bc^3$$

$$(8) \quad a^{-(3+j)}g_j ba^{(3+j)} = c^{-(3+j)}bc^{(3+j)}, j = 1,\ldots,r.$$

As before, it is not difficult to see ([M], p.13) that $B$ is the free product with amalgamation

$$(X* <a,b_1|>)* <b_2,c|>;$$

$$V_1 = V_2,$$

where $V_1$ is the free groups generated by $b_1$ and the left hand sides of relations (5)-(8), where $b$ is replaced by $b_1$; $V_2$ is the free group generated by $b_2$ and the right hand sides of the relations (5)-(8), where $b$ is replaced by $b_2$. (We replaced $b$ by two different symbols $b_1$, $b_2$ but added the equation identifying these symbols.) It is clear from the Mayer-Vietoris sequence that the cohomological dimension of $B$ is at most two.

Proofs of many lemmae below use the following well-known theorem due to Schreier (cf. [Baum], Theorem 14 on p. 59): Let $D, E$ be two groups and $G = D *_H E$ be their free product with amalgamation. 1) Let $x \in G$ be a product $g_1 \ldots g_k$ of elements $g_i$ from either $D - H$ or $E - H$ such that for any $i$ $g_i$ and $g_{i+1}$ are not both in $D - H$ or $E - H$. Then $x \neq e$ in $G$; 2) $D$ and $E$ embed into $G$; 3) $D \bigcap E = H$.

**Lemma 5.1.** Let $x \in B$ be any element which does not belong to the infinite cyclic subgroup generated by $c$. Then $x$ and $c$ are independent generators of a free subgroup of $B$.

**Proof:** Note that $c \in <b_2,c|>$ is not in $V_2$. The mentioned Schreier theorem implies this result if $x \in B - V_1$. If $x \in V_1 = V_2 \subset <b_2,c|> \subset B$, then the result follows from the fact that $<b_2,c|>$ is a free group.

**Corollary 5.1.1.** For any non-zero $m$ $c^m \in B$ is not the square of any element of $B$ which is not a power of $c$.

**Proof:** Obvious.

Let $G = <A,B|g = c>$. Let $H$ be the infinite cyclic subgroup of $G$ generated by the image of $g$ (or $c$).

**Corollary 5.1.2.** $c \in G$ is not the square of any element of $G$.

**Proof.** We know from our assumptions about $A$ and Corollary 5.1.1 that $c \in G$ is not the square of any element of torsion-free subgroups of $G$ $A$ and $B$. Let $c = x^2$, where $x \in G - (A \bigcup B)$. Then $x$ can be represented as a product of at least two elements $g_i$ from $A - H$ and $B - H$ such that for any $i$ $g_i$ and $g_{i+1}$ do not belong both to either $A - H$ or $B - H$. Now the application of the above-mentioned Schreier theorem leads to a contradiction.

**Corollary 5.1.3.** For any integer $m$ $c^m \in G$ is not the square of any element of $G - A$ conjugate to an element of $A$ or of $B$.

**Proof.** Corollary 5.1.3 easily follows from the Schreier theorem and Lemma 5.1.

**Lemma 5.2.** Let $v \in G$ be any element which does not belong to $A$. Then for any non-zero integer $m$ $v$ and $g^m \in G$ generate a free subgroup of $G$.

**Proof.** Easily follows from Lemma 5.1 and the mentioned Schreier theorem.

**Corollary 5.2.1** Any non-cyclic virtually nilpotent subgroup of $G$ is conjugate to a virtually nilpotent subgroup of $A$ or a virtually nilpotent subgroup of $B$.

**Proof:** Let $G_1$ be a virtually nilpotent non-cyclic subgroup of $G$. Then Theorems 6,7 of [KS] imply that $G_1$ is one of the following: (i) a subgroup of a conjugate of $A$; (ii) a subgroup of a conjugate of $B$; (iii) a subgroup conjugate to a subgroup generated by a power of $g$ and by some element not in $A$ or $B$; (iv) $< H_2 * H_3 | H_4 >$, where $H_2$ and $H_3$ are subgroups of conjugates of $A$ or $B$, $H_4$ is a subgroup of a conjugate of $H$, and $H_4$ is a subgroup of index two in both $H_2$ and $H_3$. However, Lemma 5.2 implies that the case (iii) is impossible, and Corollary 5.1.3 implies that in the case (iv) $< H_2 * H_3 | H_4 >$ is a subgroup of a conjugate of $A$. □

**Remark.** Instead of Theorems 6,7 of [KS] one can use here also a simpler and better known classical result by H. Neumann on subgroups of amalgamated free products (cf. [LS], Theorem 6.6 of Chapter IV). This result immediately implies that a conjugate of $G_1$ must have a non-trivial intersection with $H$ if $G_1$ is non-cyclic and virtually nilpotent and (i), (ii) do not hold. Now Lemma 5.2 can be used to demonstrate that $G_1$ contains a subgroup as in (iii). But then $G_1$ has exponential growth and, therefore, is not virtually nilpotent. This contradiction proves Corollary 5.2.1.

**Lemma 5.3.** The intersection of any two different conjugates of $A$ is the trivial element.

**Proof.** It is sufficient to prove this Lemma in the case when one of the conjugates of $A$ is $A$ itself. Assume that a non-trivial element $a \in A$ belongs to the intersection of $A$ and its conjugate by an element of $G - A$. Then $a = \prod_{i=1}^{k} g_i$, where each $g_i$ belongs either to $A - H$ or to $B - H$, at least one $g_i \in B - H$, and for any $i$ $g_i$ and $g_{i+1}$ are not both in either $A - H$ or $B - H$. This immediately leads to a contradiction with the mentioned Schreier theorem. □

**Lemma 5.4.** The intersection of any two conjugates of $B$ is the trivial element.

**Proof.** Similar to the proof of Lemma 5.3.

**Lemma 5.5.** An intersection of a conjugate of $A$ and a conjugate of $B$ is conjugate to a subgroup of $H$.

**Proof.** Follows from the Schreier theorem (Theorem 14 on p.59 of [Baum]) similarly to the proof of Lemma 5.3.

**Lemma 5.6.** The normalizer of any subgroup of $H$ in $G$ coincides with the normalizer in $A$ of the same subgroup regarded as a subgroup of $A$.

**Proof.** The proof easily follows from the Schreier theorem and Lemma 5.1.

**Lemma 5.7.** (a) The normalizer in $G$ of a subgroup of $A$ not contained in $H$ is a subgroup of $A$.

(b) The same is true for $B$.

**Proof.** Follows easily from the Schreier theorem.

Following [Gr 2], Section 6, denote for any group $T$ the complex of virtually nilpotent subgroups of $G$, defined as in section 6 of [Gr 2], by $Q_T$. Also, let $Q_T//T$ and $p : Q_T//T \longrightarrow K(T,1)$ be defined as in section 6, p.79, of [Gr 2].

**Lemma 5.8.** Let $A$ and $B$ be as above. Let $h \in H_n(A)$, $(n \geq 3)$, be a non-trivial homology class, and $H \in H_n(G)$ be the image of $h$ under the homomorphism induced by the canonical inclusion $A \longrightarrow G$. If $h$ is not in the image of the homomorphism $p_* : H_n(Q_A//A) \longrightarrow H_n(A)$, then $H$ is not in the image of the homomorphism $p_* : H_n(Q_G//G) \longrightarrow H_n(G)$.

**Proof:** First, recall the definitions of $Q_G, Q_G//G$ and $p$ from [Gr 2], Section 6. For any group $K$ the complex $Q_K$ is defined as follows: Its zero dimensional simplices correspond to non-trivial elements of $K$. For any finite set of non-trivial elements of $G$ generating a virtually nilpotent subgroup of $K$ $Q_K$ contains a simplex with vertices corresponding to these elements. $K$ acts on $Q_K$ by conjugations. Consider the product of $Q_K$ and the universal covering $EK$ of $BK$. Consider the diagonal action of $K$ on this product $Q_K \times EK$. The quotient of this product with respect to this diagonal action is, by definition, the homotopy quotient $Q_K//K$. The map $p$ of $Q_K//K$ is induced by the projection of $Q_K \times EK$ on the second factor. The definition of $K$ implies that we can think of $Q_K//K$ as glued from $BN_i$, where $N_i$ runs over normalizers of all virtually nilpotent subgroups of $K$. (We have one normalizer for any conjugate class of virtually nilpotent subgroups of $K$.) The intersection of $BN_i$ and $BN_j$ is either empty or is contained in the classifying space of the normalizer of the intersection of virtually nilpotent subgroups $K_i$ and $K_j$ of $K$ such that $N_i$ is the normalizer of $K_i$ and $N_j$ is the normalizer of $K_j$.

Having this in mind, we can easily see what is the geometry of $Q_G$ and $Q_G//G$ in the considered case. Corollary 5.2.1 and Lemmae 5.3,5.4,5.5 imply that any connected component of $Q_G$ is either made of copies of $Q_A$ and $Q_B$ in such a way that different copies of $Q_A$ do not intersect, different copies of $Q_B$ do not intersect, and copies of $Q_A$ and $Q_B$ can intersect only along a copy of $Q_H$ (if they intersect at all) or corresponds to an infinite cyclic subgroup of $G$. Note that the normalizer of any subgroup of $B$ which is not a subgroup of $H$ is a subgroup of $B$. Similarly, the normalizer of any subgroup of $A$ is a subgroup of $A$. Therefore, $Q_G//G$ will be union of circles, corresponding to the conjugacy classes of cyclic subgroups of $G$ not conjugate to a subgroup of $A$ or $B$, and also of $Q_A//A$ and $Q_B//B$, where the intersection of $Q_A//A$ and $Q_B//B$ is made of the classifying spaces of normalizers of subgroups of $H$ in $A$ (we are using Lemma 5.6 here). However, Lemma 5.1 implies that the intersection of $H \subset B$ and any virtually nilpotent subgroup of $B$ which is not a subgroup of $H$ is the trivial element. Thus, the part of $Q_B//B$ corresponding to subgroups of $H$ is a union of connected components of $Q_B//B$. Denote the union of remaining connected components of $Q_B//B$ by $(Q_B//B)_0$. Therefore $Q_G//G$ is the union of $Q_A//A$ and $(Q_B//B)_0$, which are pairwise disjoint, and circles, corresponding to the conjugacy classes of cyclic subgroups of $G$ not conjugate to a subgroup of $A$ or $B$. It remains to notice that $p_*$ for $G$ coincides with $p_*$ for $A$ on $H_n(Q_A//A)$ and $p_*$ for $B$ on

$H_n((Q_B//B)_0)$. □

Now we can finish the proof of (v.ii). Let $g_1, \ldots g_l$ be the system of generators of $A_0$ such that for any $i$ $g_i$ is an element of infinite order and not a square. For $j = 1, \ldots, l$ define recursively $A_j = < A_{j-1} * W_i^{(j)} | g_j = c^{(j)} >$. If $W_i$ is non-trivial then Corollary 5.1.2 and Lemma 5.8 imply by induction that for any $j$ the image of $h_n$ in $H_n(A_j)$ under the homomorphism induced by the canonical inclusion $A_0 \longrightarrow A_j$ is not in the image of the homomorphism $p_* : H_n(Q_{A_j}//A_j) \longrightarrow H_n(A_j)$. (On the $j$th step we apply Lemma 5.8 to $A = A_{j-1}, B = W_i^{(j)}$. $g$ is the image of $g_j$ in $A_{j-1}$. Corollary 5.1.2 implies that $g_j$ is not a square in $A_{j-1}$. The properties (5'.1) and (5'.2) of $A_0$ are used as the base of induction.) It remains to observe that $L_i = A_l$.

**Step 2. Construction of $S_i$.** Now we are going to explain an algorithm which constructs a sequence of smooth $n$-dimensional homology spheres $S_i$ such that $\pi_1(S_i) = G_i$, and the image of the fundamental homology class of $G_i$ under the homomorphism induced by the classifying map is $\bar{h}_i \in H_n(BG_i)$. First, observe that if $D$ is any perfect group and $C$ its perfect subgroup, then the universal central extension $\bar{C}$ of $C$ can be mapped into the universal central extension $\bar{D}$ of $D$ by a homomorphism $\phi$ such that the resulting diagram of groups $\bar{C}$, $C$, $\bar{D}$, $D$ becomes commutative. Indeed, assume that $C$ has generators $f_i$, $i = 1 \ldots, k$ and relators $r_j$, $j = 1, \ldots m$. Without any loss of generality we can assume that these generators and relators are among generators and relators of $D$. Now $\bar{C}$ has the same generators as $C$ and it sets of relators consists of all commutators $[r_i, f_j]$ between relators and commutators of $C$ and relators $\lambda_i$ such that $\lambda_i$ regarded as an element of the free group $F$ generated by the generators $f_i$ is a product of relators of $C$ and $(\lambda_i)^{-1} f_i$ is the product of commutators of elements of $F$. It is clear that we can ensure that these relators are among the relators of $\bar{D}$. Now we see that the homomorphism which maps the generators $f_i$ of $\bar{C}$ into the same generators regarded as generators of $\bar{D}$ is the required homomorphism $\phi$. It is easy to see looking at the exact sequences (***) for $C$ and $D$ that if the homomorphism $H_2(C) \longrightarrow H_2(D)$ induced by the inclusion is injective, then $\phi$ is an injection. Applying these observations to $C = A_0$ and $D = L_i$ we obtain an explicit homomorphism $f$ of $\bar{A}_0$ into $G_i$ which is injective for any $i$ such that $L_i$ is non-trivial.

Now recall that $B\bar{A}_0^+ \longrightarrow BA_0^+$ is the fibration with the fiber $BH_2(A_0)$ (cf. [Ros 2], Remark 5.2.9). (Here $\bar{A}_0$ denotes the universal central extension of $A_0$.) Therefore for the considered values of $n$ $\pi_n(B\bar{A}_0^+)$ is isomorphic to $\pi_n(BA_0^+)$. Now it follows from the commutative square formed by $\pi_n(B\bar{A}_0^+)$, $H_n(B\bar{A}_0^+)$, $\pi_n(BA_0^+)$, $H_n(BA_0^+)$ that $h_n^+$ is the image of some spherical homology class $s_n^+ \in H_n(B\bar{A}_0^+)$. (The horizontal arrows in this diagram are induced by the homomorphism $\bar{A}_0 \longrightarrow A_0$, and vertical arrows are induced by the Hurewicz homomorphisms.) Therefore $h_n$ is the image of the homology class $s_n \in H_n(B\bar{A}_0)$ under the homomorphism induced by the central extension $\bar{A}_0 \longrightarrow A_0$. (Observe that the Mapping Theorem on p. 40 of [Gr 1] implies that $s_n$ and $h_n$ have the same simplicial norm.)

We start our construction from taking a fixed smooth homology sphere $\Sigma^n$ such that its fundamental group is the universal central extension of $A_0$ and the image of its fundamental homology class is the just constructed class $s_n \in H_n(B\bar{A}_0)$. (The existence of some such homology sphere follows from the results of Hausmann and Vogel ([Haus]) cited above.) Since $\Sigma^n$ does not depend on $i$ we can regard some finite description of $\Sigma^n$ as a part of the algorithm. For any $i$ we construct a smooth homology sphere $\Sigma_i^n$ with the fundamental group $G_i$ embedded as a hypersurface in $\mathbb{R}^{n+1}$ such that the image of its fundamental homology class under the homomorphism induced by the classifying map into $BG_i$ is zero. (The construction of $\Sigma_i^n$ starts from the Dehn construction of a smooth hypersurface $T_i$ in $\mathbb{R}^{n+1}$ with the fundamental group $G_i$ and such that $H_3(T_i) = \ldots = H_{n-3}(T_i) = 0$. Then one realizes all generators of the second homology group of this manifold by embedded 2-spheres and kills them by the obvious surgeries. A detailed description of this construction can be found in the Appendix of [N1]. The image of the fundamental homology class of $\Sigma_i^n$ in $H_n(BG_i)$ is zero because of the following reason: This manifold can be described as the smoothed out boundary of a tubular neighborhood of an embedded in $\mathbb{R}^{n+1}$ acyclic 3-dimensional complex $K$ with the fundamental group $G_i$. The map which sends any point of the manifold to the closest point of $K$ induces the isomorphisms of the fundamental groups and sends the fundamental homology class to zero. Now the composition of this map with the obvious embedding of $K$ into $BG_i = K(G_i, 1)$ will be (homotopy equivalent to ) the classifying map for $\Sigma_i^n$.) Construct the connected sum of $\Sigma^n$ and $\Sigma_i^n$. The result will be a smooth homology sphere $D_i$ with the fundamental group $\bar{A}_0 * G_i$. Kill by the obvious surgeries the elements of this group of the form $(x_j)^{-1} f(x_j)$, where $x_j$ runs over all generators of $\bar{A}_0$ and $f$ is the homomorphism of $\bar{A}_0$ into $G_i$ introduced above. (It sends the generators $\bar{A}_0$ to the generators of $G_i$ with the same names. The generators of $\bar{A}_0$ have the same names as generators of $A_0$. Since $A_0$ does not depend on $i$, we can regard some its finite presentation as a set of data known to the algorithm. If $G_i$ is non-trivial, then $f$ is an inclusion.) Obviously the fundamental group of the resulting manifold $S_i^{n*}$ will be $G_i$. Now we are going to describe a system of 2-cycles generating $H_2(S_i^{n*})$. Denote 2-cells used to kill $(x_j)^{-1} f(x_j)$ by $d_j$. Denote 2-cells corresponding to the relators $\lambda_j$ of $\bar{A}_0$ and $G_i = \bar{L}_i$ in $\Sigma_i^n$ by $c_{2j}$ and in $\Sigma^n$ by $c_{1j}$. (Recall, that $\lambda_j$ is the word which belongs to the subgroup generated by relations of $A_0$ in the free group $F$ generated by the generators of $A_0$ and $x_j(\lambda_j)^{-1} \in [F, F] \subset F$. The cells $c_{2j}$ are constructed during the construction of $\Sigma_i^n$. On the other hand the specific choice of $c_{1j}$ is not important. We can assume that the algorithm "knows" a 2-complex in $\Sigma^n$, where all 1-cells correspond to generators of $\bar{A}_0$ and 2-cells correspond to relators of $\bar{A}_0$, and all these cells are marked by the names of the corresponding generators and relators.) Now it is easy to see that the chains $c_{1j} + d_j - c_{2j}$ are cycles generating $H_2(S_i^{n*})$. Since $\pi_1(S_i^{n*}) = G_i$ has the trivial second homology group the Hopf theorem implies that all two-dimensional homology classes of $S_i^{n*}$ are spherical. Obviously, they can be realized by embedded 2-dimensional spheres. Knowing this fact a priori we can find some such system of embedded spheres corresponding to the generators and kill them by surgeries. As the result we obtain a smooth homology sphere $S_i$ with the desired properties. (In the case, when $n = 5$ we must check that we are not going to get any new 2-cycles as the result of these surgeries. This would follow from the existence for any of the 2-cycles $\sigma_j$ in $S_i^{n*}$ we are going to kill of a 3-cycle $\tau_j$ representing a non-trivial element in $H_3(S_i^{n*})$ such that the intersection number of $\sigma_j$ and $\tau_j$ is one. (Or, more formally, the cup product of cohomology classes dual to the homology classes represented by $\sigma_j$ and $\tau_j$ is the fundamental cohomology class.) To construct $\tau_j$ note that $\sigma_j$ is formed by a 2-chain $c_D$ in $D_i$ filling an embedded circle $C_j$ representing $(x_j)^{-1}f(x_j)$ and a new 2-dimensional chain $c_S$ bounded by the same circle which is the axis of the handle added to kill the element $(x_j)^{-1}f(x_j)$. Consider a small 3-dimensional sphere around $C_j$ in $D_i$. Assume here that the radius of the tubular neighborhood of $C_j$ deleted during the surgery killing $x_j^{-1}f(x_j)$ is less than the distance from $C_j$ to the considered 3-sphere. It is easy to see that the 3-dimensional sphere survives this surgery and it is easy to see that it intersects $c_D$ at one point and does not intersect $c_S$. Therefore we can define $\tau_j$ as the 3-dimensional cycle in $S_i^{n*}$, corresponding to this three-dimensional sphere.)

**Step 3. From the algorithmic unsolvability of the homeomorphism problem to the algorithmic unsolvability of the diffeomorphism problem.** Note that
1) If $G_i$ is trivial, then $S_i$ is homeomorphic to $S^n$;
2) If $G_i$ is non-trivial, then $S_i$ does not admit a Riemannian metric of non-negative scalar curvature, and for any Riemannian metric on $S_i$ of Ricci curvature $\geq -(n-1)$ the volume is not less than a positive constant not depending on $i$. (Below we will call this property of $S_i$ "pseudohyperbolicity".)
These properties and the properties of groups $G_i$ described above imply that there is no algorithm which decides whether or not a given smooth homology sphere which is a priori known to be either a homotopy sphere or a pseudohyperbolic manifold is a homotopy sphere. We must now get rid of the word "homotopy" here. Here are two different ways to do that. One way is to show that the homeomorphism to $S^n$ problem is reducible to the diffeomorphism problem. Indeed, form connected sums of $S_i$ with all distinct smooth $n$-dimensional homotopy spheres. (For every fixed $n \geq 5$ there are only finitely many distinct smooth homotopy spheres. We can regard their descriptions as data known to the algorithm.) Now $S_i$ is homeomorphic to $S^n$ if and only if at least one of these connected sums is diffeomorphic to $S^n$. If $S_i$ is not homeomorphic to $S^n$ then all connected sums of $S^n$ and smooth $n$-dimensional homotopy spheres will be pseudohyperbolic by the same reason as $S_i$. Now let us apply in parallel the algorithm solving the diffeomorphism problem to all these connected sums. If $S_i$ is not homeomorphic to $S^n$, then the answer in all cases will be negative. If $S_i$ is homeomorphic to $S^n$, then in some cases the data for the algorithm will be inappropriate (that is, it will get a manifold which is not pseudohyperbolic but is not diffeomorphic to $S^n$) and there is no way to say what will be the outcome of the work of the algorithm (for example, it might work the infinite time) but for at least one connected sum the answer will be "yes", and this answer will be found in a finite time. Thus, we need just to wait for either the first positive answer or till the algorithm will produce the negative answers for all connected sums.

Another way to proceed is to introduce a small change in the construction of $S_i$ to ensure that $S_i$ is homeomorphic to $S^n$ if and only if $S_i$ is diffeomorphic to $S^n$. It is well-known (and not difficult to see) that for any smooth $n$-dimensional homology sphere $X$ there exists a smooth $n$-dimensional homotopy sphere $Y$ such that the connected sum of $X$ and $Y$ is embeddable into $\mathbb{R}^{n+1}$. Thus, taking if necessary a connected sum with an appropriate homotopy sphere we can assume that the homology sphere $\Sigma^n$ introduced at the beginning of Step 2 and used in the construction of $S_i$ is embeddable into $\mathbb{R}^{n+1}$. Since $\Sigma^n$ is the same for all $S_i$ we can assume that the algorithm constructing $S_i$ "knows" some embedding of $\Sigma^n$ into $\mathbb{R}^{n+1}$. Since all homology spheres $\Sigma^n_i$ introduced on Step 2 were constructed as hypersurfaces of $\mathbb{R}^{n+1}$, the connected sums of $\Sigma^n$ and $\Sigma^n_i$ are hypersurfaces in $\mathbb{R}^{n+1}$. It remains to check that all surgeries performed at the end of Step 2 on the connected sums of $\Sigma^n$ and $\Sigma^n_i$ can be made inside $\mathbb{R}^{n+1}$. This is obvious if $n \neq 5$ since we added only 2-handles and 3-handles. If $n = 5$ then we can still add 2-handles without any problems, but we need the Whitney trick and some care to add the 3-handles. (We refer the reader to p. 89 of the Appendix of [N1] where it is explained in more details how to perform surgeries on a hypersurface inside $\mathbb{R}^6$ in a similar situation.) In the rest of the proof we will denote by $S_i$, $i = 0, 1, \ldots$ the elements of the just constructed sequence of connected sums and not of the "original" sequence of homology spheres constructed at the end of Step 2.

**Step 4. Complexity of the constructed homology spheres.** Now we would like to choose a Riemannian metric $\mu_i$ on the constructed smooth homology spheres $S_i$ such that the convexity radius is bounded from below by 1, the absolute values of sectional curvature are bounded by one, the diameter does not exceed $C_1(\ln(i+1))^{\frac{1}{n}}$ and the volume is between $C_2 \ln(i+1)$ and $C_3 \ln(i+1)$ for some positive $C_1$, $C_2$, $C_3$ not depending on $i$. The constructed homology spheres $S_i$ were built on Step 2 by surgeries on the connected sum of homology spheres $\Sigma^n$ and $\Sigma^n_i$. Note that the complexity of any fixed Riemannian metric on $\Sigma^n$ does not depend on $i$ at all. Also note that an almost obvious smooth version of the construction of PL-analogues of homology spheres $\Sigma^n_i$ in the section 2 of [N2] yields Riemannian metrics on $\Sigma^n_i$ with the desired properties. Take a Riemannian connected sum of $\Sigma^n_i$ endowed with such a Riemannian metric and $\Sigma^n$ with an arbitrary metric (the same for all $i$). We will not have any problems with choosing Riemannian metrics on the added handles and obtaining after an appropriate rescaling the desired Riemannian metric on $S_i$ if we are able to prove that the surgeries on $\Sigma^n \# \Sigma^n_i$ during the construction of $S_i$ on Step 2 could be performed in a fashion independent of $i$. To achieve this goal note that the same construction as in the proof of Theorem 2.1 in [N2] implies that the required surgeries correspond to any representation of elements $e_j = \lambda_j(x_j)^{-1} f(x_j)(f(\lambda_j)^{-1})$ in the free group $F$ generated by all generators of $\pi_1(\Sigma^n)$ and $\pi_1(\Sigma^n_i)$ as the product of commutators of the relators of $< \pi_1(\Sigma^n) * \pi(\Sigma^n_i) | (x_1)^{-1} f(x_1), \ldots, (x_m)^{-1} f(x_m) >$ and some elements of $F$, where $x_1, \ldots, x_m$ is a basis of $A_0$. Note that the set of elements $e_j$ does not depend on $i$. If we will show that these elements can be represented as the products of commutators of relators and arbitrary words in a way not depending on $i$ we will be done. Here is a specific way to represent these elements as the required products: Recall that according to the definition of $\lambda_j$ $x_j = c_j \lambda_j$, where $c_j = \prod_{l=1}^{K_j} [g_{lj1}, g_{lj2}]$. (Here $g_{ljq}$ are arbitrary words in $x_1, \ldots x_m$ and their inverses.) Thus, $e_j = (\prod_{l=1}^{K_j}[g_{lj1}, g_{lj2}])^{-1} \prod_{l=1}^{K_j} [f(g_{lj1}), f(g_{lj2})]$. Note that for any generator $x_k$ $z_k = (x_k)^{-1} f(x_k)$ is a relator in the considered finite presentation, and $f(x_k) = x_k z_k$. We can use these identities to replace every expresssion $f(g_{ljq})$ in the last formula for $e_j$ by the product of powers of generators $x_k$ and relators $z_k$. Move all powers of relators $z_k$ to the right side of the formula for $e_j$. A lot of commutators of relators $z_k$ and words in the generators $x_1, \ldots x_m$ and their inverses will be formed in the process. At the left side of the formula we obtain the product $(\prod_{l=1}^{K_j}[g_{lj1}, g_{lj2}])^{-1} \prod_{l=1}^{K_j} [g_{lj1}, g_{lj2}] = e$. At the right side of the formula we will get a product of powers of relators $z_k$. Since every $f(g_{ljq})$ enters the formula for $e_j$ the same amount of times as $f(g_{ljq})^{-1}$ the sum of exponents of $z_k$ will be equal to 0 for any $k$. Thus, permuting $z_1, \ldots$ we can cancel all of them creating in the process several new commutators of relators $z_1, \ldots$ with some words. As the result we obtain the identity between $e_j$ and a product of commutators of the desired structure in the free group generated by the generators in the considered finite presentation. It is clear that nothing in this identity depends on $i$. This completes the Step 4.

**Step 5. Proof of the statement about the time-bounded Kolmogorov complexity.** According to the Barzdin theorem (cf. [B], [ZL], Theorem 2.5, [LV], [N2]) there exists a Turing machine such that for any recursive time-bound $t$ the time-bounded Kolmogorov complexity $K^{(t)}(Halt(T_0), N)$ of the halting problem for $T_0$ is bounded from below by $2^N/c(t) - const$, where $N$ is the length of inputs for $T_0$, $c(t)$ is a positive constant depending only on $t$. Above we described an effective reduction of the halting problem for an arbitrary Turing machine $T$ and inputs of length $\leq N$ to the algorithmic problem which of the homology spheres $S_1, \ldots S_I$ is diffeomorphic to $S^n$ for $I \leq Const 2^N$, where $Const$ is equal either to the number of smooth $n$-dimensional homotopy spheres or to one depending on which of two ways to proceed one chooses on Step 3. (Recall that the construction of the sequence $S_i$ started in Step 1 from the choice of an arbitrary Turing machine $T$ with unsolvable word problem.) Now let this Turing machine be $T_0$. Such choice of the Turing machine immediately implies the statement about the time-bounded Kolmogorov complexity in the text of Theorem 1 for the resulting sequence $\{S_i\}$.

**Step 6. From $S^n$ to arbitrary compact $n$-dimensional manifolds.** In order to prove Theorem 1.C we can just replace the constructed sequence $\{(S_i, \mu_i)\}$ by the sequence of appropriately rescaled Riemannian connected sums of $(S_i, \mu_i)$ with $M_0^n$ endowed with an arbitrary Riemannian metric. The only statement which is not obvious is the statement that in the case when $M_0^n$ is spin and $S_i$ is nonsimply-connected the connected sum of $M_0^n$ and $S_i$ does not admit Riemannian metrics of non-positive scalar curvature even in the case when the fundamental group of $M_0^n$ is not known to satisfy the Strong Novikov Conjecture (SNC 2 in terminology of [Ros 1]). This statement can be proven by contradiction. Assume that the connected sum admits a Riemannian metric of non-negative scalar curvature. Kill the fundamental group of $M_0^n$ by surgeries (but preserving the fundamental group of $S_i$). The resulting manifold still admits a Riemannian metric of non-negative scalar curvature ([GL 1]). But now the same theorem in [Ros 1] which was previously used to prove that $S_i$ does not admit a Riemannian metric of non-negative scalar curvature implies that the resulting manifold does not admit such a Riemannian metric. This fact provides the desired contradiction. □

### 3. Some useful discrete groups.

Here we are going to construct the groups $A_0$ with the properties postulated in the Proposition 4. We are going to present three essentially different constructions. The quickest way to prove Proposition 4 is to use the second construction below. This construction has an extra advantage that it yields $A_0$ which satisfies all conditions of Proposition 4 (that is, (1)-(5), (5'.1), (5'.2)) and is obtained as the amalgamated free product of the fundamental group of a compact Kähler manifold of negative sectional curvature with several copies of an explicit acyclic group. The main idea of the second construction is to use a certain class of arithmetic groups investigated in the Clozel paper [Cl]. These groups are fundamental groups of compact manifolds of negative sectional curvature and, therefore, all their non-trivial homology classes of infinite order in dimensions $> 1$ satify (5). On the other hand these arithmetic groups satisfy particularly strong vanishing theorems. As a corollary, killing their first homology groups by passing to an amalgamated free product with acyclic groups where all amalgamated subgroups are cyclic we obtain groups such that for an appropriate homology class (3) can be proven just by an application of the Serre generalization of the Hurewicz theorem. The Clozel paper [Cl] is based on an impressive amount of previous results obtained by different authors and related to the Langlands program. Therefore we decided to look also for simpler constructions using less machinery.

The third construction also uses some number theory and representation theory but of a less recent vintage. We start from the most well-known fundamental group of a compact orientable $n$-dimensional manifold of constant negative sectional curvature. This group is an appropriate torsion-free subgroup of finite index of $SO_F(J)$, where $J$ is the ring of integers of $\mathbb{Q}(\sqrt{2})$, and $F = x_1^2 + \ldots + x_n^2 - \sqrt{2}x_{n+1}^2$. We modify the construction of this group by replacing $J$ by $\mathbb{Z}[\sqrt{2}, \frac{1}{5}, \frac{1}{13}]$ (or even by $\mathbb{Z}[\sqrt{2}, \frac{1}{5}]$). The resulting S-arithmetic group will still be formed by hyperbolic isometries of the hyperbolic space $H^n$. This fact can be used to prove (5'.1) for an appropriate $n$-dimensional homology class dual to a nontrivial cohomology class coming from the $n$-th continuous cohomology group of $SO(n,1)$ with real coefficients. The main idea of the construction is in the observation that the results of Chapter XIII of [BorWal] imply vanishing of cohomology groups with real coefficients of the constructed S-arithmetic groups in low dimensions which enables us to apply a generalization of Hurewicz theorem to prove the property (3) in Proposition 4. It remains to get rid of the possibly nontrivial (finite) first homology group of the constructed $S$-arithmetic group. This is done by forming an amalgamated free product of the S-arithmetic group with several copies of an explicit acyclic group, where all amalgamated groups are cyclic. Now the (repeated) application of Lemma 5.8 from section 2 implies (5'.1) for the amalgamated free product.

Our first construction is the most elementary. Unfortunately it only works for even $n$ (although we sketch a modification of this construction which might work for odd $n$ as well). For even $n$ this construction yields pairs $A_0, h_n \in H_n(A_0)$ satisfying conditions (1)-(5) of the Proposition 4. Here are the main ideas of the first construction for even $n$. To prove (5) we use the theorem on p. 23 of [Gr 1] asserting that real cohomology classes of the classifying space of an algebraic group $G \subset GL_N(\mathbb{R})$ defined over $\mathbb{R}$ and endowed with the discrete topology are bounded if they are in the image of $H^*(BG)$ under the homomorphism induced by the inclusion $G^{discrete} \longrightarrow G$. Moreover, such classes belong to the continuous cohomology with real coefficients and, thus if they are nontrivial, then they are mapped to nontrivial cohomology classes under the homomorphism into the cohomology of a discrete cocompact subgroup $\Gamma$ of $G$ ([Mil 2], [Bor 2]). In fact, this result holds in an appropriate range of dimensions even in the case when $\Gamma$ is a not necessarily cocompact lattice. Combining the mentioned result on p. 23 of [Gr 1] with results stated in the Appendix of [Mil 2] and with calculations of real cohomology of specific arithmetic groups in [Bor 2] we are able to conclude that the groups $\Gamma_N$ defined as $Sp_{N,N}(\mathbb{Z})$ (for $n$ divisible by four) and $Sp_{2N}(\mathbb{Z})$ (for even $n$ not divisible by four) have a non-trivial $n$-dimensional homology class satisfying (5) providing that $N$ is sufficiently large (e.g. $N \geq 4(n+1)$).) In order to kill the first homology group of $\Gamma_N$ and obtain a finitely presented group satisfying (2) we can take an amalgamated free product of $\Gamma_N$ with certain groups which are either acyclic or have nontrivial homology only in dimension two, and where all amalgamated subgroups are cyclic. In order to satisfy (3) we notice that $BSL(\mathbb{Z})^+$ and the results of the application of the +-construction to classifying spaces of direct limits $\varinjlim \Gamma_i$ of many similar sequences of arithmetic groups $\Gamma_i$ are H-spaces (cf. [Ros 2], [Bor 2]). As a corollary their cohomology rings (with rational coefficients) are generated by indecomposable generators dual to rational homology classes in the image of the Hurewicz homomorphism from tensor products of homotopy groups with $\mathbb{Q}$ (cf. [Bor 2]. Indeed, the cohomology algebra with real coefficients of an H-space is a tensor product of a polynomial algebra generated by some elements of even degrees and an exterior algebra generated by some elements of odd degrees. One can realize any of these indecomposable generators $\{x_i\}$ by maps from the H-space into corresponding Eilenberg-Maclane space $K(\mathbb{Q}, \deg x_i)$. Combine these maps into one map from the H-space into the product of all these Eilenberg-Maclane spaces. This map induces the isomorphism of all rational cohomology groups and as a corollary it is a rational homotopy equivalence.) Now it is natural to try to use stabilization theorems to prove that the indicated above $n$-dimensional homology classes of $B\Gamma_N^+$ are spherical. We simplify our calculations and avoid some minor technical problems by introducing and using a "rational +-construction" denoted below by the superscript $R$. Roughly speaking, it is a composition of the +-construction and the rational localization but unlike the +-construction this construction is always defined for CW-complexes with a finite first homology group.

**First construction.** First, we are going to introduce a useful functor from the homotopy category of CW-complexes with finite (or torsion) first homology group into the homotopy category of simply connected rational spaces. We will call it "rational +-construction" and denote the result of its application by the superscript $R$. It can be defined as follows: For a CW-complex $X$ with a finite first homology group we first kill its fundamental group by adding new 2-dimensional cells (one cell for any generator). The second homology group with rational coefficients of the resulting space $\tilde{X}$ is the direct sum of $H_2(X, \mathbb{Q})$ and a group freely generated by the following cycles: Assume that a newly added 2-cell $C$ kills the generator $g$ of $\pi_1(X)$ such that $k[g] = e$ for some $k \neq 0$, where $[g]$ denotes the image of $g$ under the Hurewicz homomorphism $\pi_1(X) \longrightarrow H_1(X)$. Let $D$ be a 2-chain in $X$ such that its boundary is $k[g]$. Then $kC - D$ is a generator of $H_2(\tilde{X}, \mathbb{Q})$, and the Hurewicz theorem applied to the simply connected space $\tilde{X}$ implies that this generator is representable by a map of the two-dimensional sphere into $\tilde{X}$.

Let us kill now these new generators of the second homology group by adding 3-cells. As the result we will obtain a simply connected CW-complex $\check{X}$ such that $H_*(\check{X}, \mathbb{Q}) = H_*(X, \mathbb{Q})$. Define $X^R$ as the rational localization of $\check{X}$. (The construction and the properties of rational localization can be found, for example, in [GM], p. 90, where it is called localization at 0 or in [FHT]). Here note only that for any $i$ $\pi_i(X^R) = \pi_i(\check{X}) \otimes \mathbb{Q}$ and $X^R$ comes with a canonical continuous map $\check{X} \longrightarrow X^R$ inducing the isomorphism of all homotopy groups tensored with $\mathbb{Q}$.) Clearly $X$ embeds into $\check{X}$. Therefore there is a canonical map $i_X : X \longrightarrow X^R$ inducing the isomorphism of homology groups with rational coefficients. Also, it is clear that if $\pi_1(X)$ is perfect, then $X^R$ is simply the rational localization of $X^+$.

**Proposition 7.** Let $X$ be a CW-complex such that $H_1(X)$ is finite. Then

(1) For any simply connected CW-complex $Z$ such that all homotopy groups of $Z$ are $\mathbb{Q}$-vector spaces any continuous map $\phi : X \longrightarrow Z$ can be passed through $X^R$, that is, there exists a map $\psi : X^R \longrightarrow Z$ such that $\psi \circ i_X$ is homotopic to $\phi$. Moreover, $\psi$ is unique (up to a homotopy).

(2) For any morphism $F : X \longrightarrow Y$ there exists a unique morphism $F^R : X^R \longrightarrow Y^R$ such that $i_Y \circ F$ is homotopic to $F^R \circ i_X$. (Of course, here $H_1(Y)$ is also assumed to be finite.)

**Proof.** Part (1) easily follows from the obstrucion theory. Indeed, all homotopy groups of $Z$ are rational, but $X$ and $X^R$ are rationally homology equivalent. Therefore all homology groups of the pair $(X^R, X)$ with coefficients in homotopy groups of $Z$ vanish. Part (2) easily follows from part (1). □

We will need the following two obvious lemmae. The first of them immediately follows from Proposition 7.

**Lemma 8.** If $\pi_1(X)$ has a perfect commutator subgroup and $H_1(X)$ is finite, then $(X^+)^R = X^R$. Moreover, $i_X$ is the composition of $i_{X^+}$ and the canonical inclusion $X \longrightarrow X^+$. Also, $i_{X^+} : X^+ \longrightarrow X^R$ induces the isomorphism of all homology groups with rational coefficients.

**Lemma 9.** Assume that $I : A \longrightarrow C$ is an inclusion of groups such that $H_1(A)$ and $H_1(C)$ are finite. Assume that this inclusion induces the isomorphism of $H_n(A, \mathbb{Q})$ and $H_n(C, \mathbb{Q})$. Let for some $h \in H_n(BA, \mathbb{Q})$ $h^R = i_{BA}(h) \in H_n(BA^R, \mathbb{Q})$ be the image of some element from $\pi_n(BA^R) = \pi_n(BA^R) \bigotimes \mathbb{Q}$ under the Hurewicz homomorphism. Then $I_*^R(h^R) \in H_n(BC^R, \mathbb{Q}) = H_n(BC^R)$ will be in the image of the Hurewicz homomorphism. If $C$ is perfect, then $BC^R$ is the rational localization of $BC^+$. In this case, $I_*^R(h^R)$ is the image under the isomorphism induced by rational localization of the class $\hat{h}^+ = (I_*(h))^+ \in H_n(BC^+, \mathbb{Q})$ such that its integer multiple is a non-zero spherical integer homology class.

After all these preliminaries we are going to describe the construction. It works only for all even $n$. However, afterwards we will outline a modification of this construction. We conjecture that this modification yields the finitely presented groups and their $n$-dimensional homology classes satisfying (1)-(4), (5'.1), (5'.2) for all odd $n$.

**The case of an even $n \geq 4$.**

We are going to construct groups $A_0$ satsfying the conditions (1)-(5) of Proposition 4. Let $n$ be an even number greater than or equal to 4. Define $A_*$ as $Sp_{2N}(\mathbb{Z})$ for some $N >> n$, if $n$ is not divisible by four, and as $Sp_{N,N}(\mathbb{Z})$ for some $N >> n$, if $n$ is divisible by four.

**Lemma.** There exist a bounded (in the sense of [Gr 1]) indecomposable class $h^n \in H^n(A_*, \mathbb{R})$.

**Proof.** In the proof of this lemma we consider only cohomology with rational coefficients and therefore will be omitting for brevity the coefficients. According to Milnor ([Mil 2], Appendix) for any real connected semi-simple Lie group $G$ there exists the following exact sequence of homomorphisms:

$$H^*(BG_{\mathbb{C}}) \longrightarrow H^*(BG) \longrightarrow H^*(BG^{discr}).$$

Here $G_{\mathbb{C}}$ denotes a connected complex Lie group whose Lie algebra is the complexification of the Lie algebra $g$ of $G$, $G^{discr}$ denotes $G$ regarded as a discrete group and "exact" means that the kernel of the second homomorphism is the ideal generated in the ring $H^*(BG)$ by the image of the first homomorphism. Further, the second homomorphism can be factored through $H^*(Inv_G A(G/K))$, where $K$ denotes the maximal compact subgroup of $G$ and $Inv_G A(G/K)$ denotes the De Rham complex of $G$-invariant smooth differential forms on $G/K$, and $H^*(Inv_G A(G/K))$ injects into $H^*(BG^{discr})$ ([Mil 2], p. 81-82). Let $\Gamma$ be an arithmetic subgroup of an almost simple Lie group $G$. Then according to Borel the composition of homomorphisms $H^*(Inv_G A(G/K)) \longrightarrow H^*(BG^{discr}) \longrightarrow H^*(B\Gamma)$ is an isomorphism in low dimensions (up to at least $(rank_{\mathbb{R}}(G)/4) - 1$. If $\Gamma$ is a discrete torsion-free cocompact subgroup of $G$ then this composition of homomorphisms is injective in all dimensions and an isomorphism in all dimensions $< rank_{\mathbb{R}}(G)$ ([Ser], p.125, [Bor 3]).) Further, $H^*(Inv_G A(G/K))$ is isomorphic to the cohomology of the compact dual of $G/K$, that is of $K_{\mathbb{C}}/K$, where $K_{\mathbb{C}}$ is the maximal compact subgroup of $G_{\mathbb{C}}$ ([Mil 2], Appendix). In particular, if $G = Sp_{2N}(\mathbb{R})$, then $G_{\mathbb{C}} = Sp_{2N}(\mathbb{C})$, $K_{\mathbb{C}} = USp_N$ (the notation of [Bor 2] for the compact symplectic group), $K = U_N$. $H^*(K_{\mathbb{C}}/K)$ coincides in lower dimensions with the polynomial algebra with exactly one indecomposable generator in every even degree not divisible by four. If $G = Sp_{N,N}$, then $K = USp_N \times USp_N$, $K_{\mathbb{C}} = USp_{2N}$, and $H^*(K_{\mathbb{C}}/K)$ coincides in lower dimensions with the polynomial algebra with exactly one indecomposable generator of any degree divisible by four. The cokernel of the restriction homomorphism $H^*(BG_{\mathbb{C}}) \longrightarrow H^*(BG)$ for $G = Sp_{2N}(\mathbb{R})$ coincides with the cokernel of $H^*(BUSp_N) \longrightarrow H^*(BU_N) = P(\{x_i, 2i | i = 1, \ldots N\})$. But $H^*(BUSp_N) = P(\{y_i, 4i | i = 1, \ldots N\})$. Therefore the cokernel will be generated by the indecomposable generators of even degrees not divisible by four. By the virtue of the discussion above this cokernel in lower dimensions is then mapped isomorphically on $H^*(BA_*)$ and this isomorphism factors through $H^*(BSp_{2N}(\mathbb{R})^{discr})$. Now Theorem on p.23 of [Gr 1] implies that the image of this cokernel in $H^*(BSp_{2N}(\mathbb{R})^{discr})$ lies in the image of bounded cohomology of $BSp_{2N}(\mathbb{R})^{discr}$. Therefore the image of this cokernel in $H^*(BA_*)$ under the described homomorphism in the lower dimensions is in the image of bounded cohomology of $BA_*$. So if $n \equiv 2 \pmod 4$ we can take for $h^n$ the image of the indecomposable element $x_n$ in $H^n(BA_*)$. Similarly, the cokernel of the ring homomorphism $H^*(B(Sp_{N,N})_{\mathbb{C}}) = H^*(BUSp_{2N}) \longrightarrow H^*(BSp_{N,N}) = H^*(B(USp_N \times USp_N))$ contains the polynomial subalgebra $P(\{x_i, 4i | i \leq N\})$. Proceeding as above one can find an indecomposable bounded cohomology class $h^n \in H^n(A_*)$, if $n$ is divisible by four, and $N >> n$. □

Denote $\varinjlim BSp_{2N}(\mathbb{Z})$ by $BSp(\mathbb{Z})$ and $\varinjlim BSp_{N,N}(\mathbb{Z})$ by $BSp_{,,,}(\mathbb{Z})$. The computations of Borel of $H^*(BSp(\mathbb{Z}))$ and $H^*(BSp_{,,,}(\mathbb{Z}))$ ([Bor 2]) imply that $h^n$ will be the image under the inclusion homomorphism of an indecomposable element $v^n$ from $H^*(BSp(\mathbb{Z}), \mathbb{Q})$ (or $H^*(BSp_{,,,}(\mathbb{Z}), \mathbb{Q})$) to $H^*(A_*)$ and that this inclusion will be an isomorphism. Consider also the dual of $v^n$ in $H_n(BSp(\mathbb{Z}))$ (correspondingly, in $H_n(BSp_{,,,}(\mathbb{Z}))$). Denote this dual by $g_n$. Now observe that the spaces $BSp(\mathbb{Z})^R$ and $BSp_{,,,}(\mathbb{Z})^R$ are H-spaces. The proof of this fact almost verbatim repeats the proof of a similar statement that spaces $BGL(\mathbb{Z})^+$ are H-spaces in [Ros 2]. The image $g_n^R$ of $g_n$ in $H_n(BSp(\mathbb{Z})^R, \mathbb{Q})$ (or $H_n(BSp_{,,,}(\mathbb{Z})^R, \mathbb{Q})$) will be dual to the indecomposable cohomology class $v^{nR}$, and therefore rational homotopy theory implies that $g_n^R$ will be in the image of the Hurewicz homomorphism. (Alternatively to prove that $g_n^R$ is in the image of the Hurewicz homomorphism from $\pi_n(BSp(\mathbb{Z})^R)$ (correspondingly, $\pi_n(BSp_{,,,}(\mathbb{Z})^R)$), we could use the known fact that the commutator groups of $Sp(\mathbb{Z})$ and $Sp_{,,,}(\mathbb{Z})$ are perfect, and therefore $BSp(\mathbb{Z})^+$ and $BSp_{,,,}(\mathbb{Z})^+$ are defined and are H-spaces. (This fact is due to Vasserstein, cf. [Bor 2], sections 11.4, 11.6, 11.7, 12.3 and references there.) The image $g_n^+$ of $g_n$ in $H_n(BSp(\mathbb{Z})^+)$ (correspondingly, in $H_n(BSp_{,,,}(\mathbb{Z})^+)$ is known to be a primitive element ([Bor 2]). Therefore it will be in the image of the Hurewicz homomorphism from the tensor product of the $n$-th homotopy group with $\mathbb{Q}$, and Lemma 8 implies the same statement for $g_n^R$.) Observe that the inclusion of $BA_*$ into $BSp(\mathbb{Z})$ (correspondingly, $BSp_{,,,}(\mathbb{Z})$) induces the isomorphism of homology groups with rational coefficients in dimensions at least up to $n$. Therefore the inclusion of $BA_*^R$ into $BSp(\mathbb{Z})^R$ (correspondingly, $BSp_{,,,}(\mathbb{Z})^R$) induces isomorphisms of all homotopy and homology groups in dimensions up to $n$. Now it is clear that the image $h_{*n}^R$ of the dual $h_{*n}$ of $h^n$ in $H_n(BA_*^R)$ is in the image of the Hurewicz homomorphism. Since $h^n$ is bounded, $h_{*n}$ has non-zero simplicial norm equal to the inverse of the norm of $h^n$ ([Gr 1]).

To complete the proof it remains to get rid of $H_1(A_*)$. Let $a_1, \ldots a_M$ be a system of elements of $A_*$ such that their images in $A_*/[A_*, A_*] = H_1(A_*)$ generate $H_1(A_*)$. Let for any $j = 1, \ldots M$ $Zj$ denotes $\mathbb{Z}$ if $a_j$ is an element of infinite order in $A_*$, and $\mathbb{Z}_{p_j}$, if $a_j$ is an element of a finite order $p_j$ in $A_*$. Denote the result of the application of the Miller "witness" construction described above to the group $Zj$ by $Y_j$. ($Y_j$ has the following finite presentation: It has four generators $a, b, c, v_j$ and relations

$$(9) \quad a^{-1}ba = c^{-1}b^{-1}cbc$$

$$(10) \quad a^{-2}b^{-1}aba^2 = c^{-2}b^{-1}cbc^2$$

$$(11) \quad a^{-3}[v_j, b]a^3 = c^{-3}bc^3$$

$$(12) \quad a^{-4}v_j ba^4 = c^{-4}bc^4,$$

If $a_j$ has a finite order $p_j$, then $Y_j$ has the fifth relation:

$$(13) \quad v_j^{p_j} = e.$$

Now let $A_0 = <A_* * Y_1 * \ldots * Y_M; a_1 = v_1, \ldots, a_M = v_M >$. It is clear that the inclusion $A_* \longrightarrow A_0$ induces the isomorphism of the $n$th homology groups and that $A_0$ is perfect.

Therefore Lemma 9 applied to $A = A_*$, $C = A_0$ implies that a non-zero integer multiple $h_n$ of the image of the generator $h_{*n}$ of $H_n(A_*)$ under the inclusion homomorphism will be spherical. Theorem on p. 55 of [Gr2] implies that the simplicial norm of $h_n$ is not less than the simplicial norm of $h_{*n}$. This completes the construction of $A_0$ for even $n$.

**The case of odd** $n \geq 5$. Let $G = SL_{\frac{n+3}{2}}(\mathbb{C})$. Alternatively, we can take $G = SL_{\frac{n+3}{2}}(\mathbb{R})$, if $n \equiv 1 \pmod{4}$, and $SO_{n+1}(\mathbb{C})$, if $n \equiv 3 \pmod{4}$. Let $A_1$ be a discrete cocompact group of $G$ constructed as in [Bor 1], $A_*$ be some its torsion-free congruence subgroup and $A$ be the acyclic group with four generators $a, b, c, v$ and four relations obtained from (9)-(12) by replacing everywhere $v_j$ by $v$. Let $A_0$ be obtained from $A_*$ by taking the amalgamated free product with several copies of the acyclic group $A$ in order to kill $H_1(A_*)$ as in the case of even $n$. Let further $h_n$ be a homology class obtained as above from the irreducible generator of the continuous cohomology $H_{ct}^n(A_*, \mathbb{R})$. (The continuous cohomology of $A_*$ is formed by harmonic $G$-invariant differential forms on $A_* \backslash G/K$, where $K$ denotes the maximal compact subgroup of $G$. Our choice of $G$ ensures that there exists exactly one (up to a multiplication by a scalar) irreducible generator of $H_{ct}^n(A_*, \mathbb{R})$ (see [Bor 2]).) We conjecture that the pair $A_0, h_n$ will satsify the conditions (1)-(4), (5'.1), (5'2) of the Proposition 4. Indeed, one can prove that $h_n^+$ is spherical by comparison with $B \xrightarrow{\lim} A_*(j)^R$ where $A_*(j)$ are formed by all $j$ exactly as $A_*$ was formed for $n$ (in particular, $A_*(n) = A_*$). (One can prove that $B \xrightarrow{\lim} A_*(j)^R$ is an H-space by the standard argument used in the algebraic K-theory to prove that $BGL(\mathbb{Z})^+$ is an H-space, cf. [Ros 2].) Therefore its real cohomology ring is generated by classes dual to real spherical homology classes. The rank of $G$ is high enough to prove using the stabilization (see [Bor 2] and [Bor 3]) and rational homotopy theory that $h_n^+$ is spherical.) The main difficulty is in the proof of (5'.1). Unfortunately, the dimension $n$ of interest for us is too low to apply Theorem 6.6 D' of [Gr 2]. After some analysis of $Q_{A_0}//A_0$ we think that the rank of $G$ is low enough to guarantee that $h_n$ is not in the image of the homomorphism $p_*$ in (5'.1), although presently we have not succeeded in verifying that.

**Second construction:** In [Cl] Clozel described a family of torsion-free discrete cocompact subgroups of $U(p,q)$ which satisfy particularly strong vanishing theorems. For example, when $p = 2n - 1$ and $q = 1$ the cohomology ring with real coefficients of any of these groups can be presented as the direct sum of the continuous cohomology which is generated by one polynomial generator in dimension 2, and another direct summand denoted $H^{var}$ which can have non-trivial elements only starting from the dimension $n$. Moreover, according to [Cl] for some groups in this class the direct summand $H^{var}$ has non-trivial elements in dimension $n$. Denote some such torsion-free discrete cocompact subgroup of $U(2n - 1, 1)$ by $\tilde{A}_0$. Since $U(2n - 1, 1)/U(2n - 1) \times U(1)$ is the Hermitian hyperbolic space, $\tilde{A}_0$ is the fundamental group of a compact manifold of negative sectional curvature and, thus, all nontrivial real homology classes of $B\tilde{A}_0$ have non-zero simplicial norms ([Gr 1], section 1.2). Now we can get rid of the (finite) group $H_1(\tilde{A}_0)$ by forming appropriate amalgamated free products with several copies of the group generated by four elements $a, b, c, v_j$ with the set of relations (9)-(12) exactly as in the first construction (recall that all amalgamated subgroups are cyclic) and then kill the second homology group by passing to the universal central extension $\bar{A}_0$ of the resulting group $A_0$. Note that by doing that we will kill all continuous cohomology but the non-trivial real cohomology class of dimension $n$ in the second direct summand will survive and will be the non-trivial real cohomology class of the lowest dimension. Consider its dual real homology class $\bar{h}_n$ of the same dimension. Now the Hurewicz theorem modulo the Serre class of finite groups implies that the image $\bar{h}_n^+$ of $\bar{h}_n$ in $H_n(B\bar{A}_0^+)$ is spherical. Note that the already mentioned Theorems on p. 55 and 40 in [Gr 1] imply that $\bar{h}_n$ will have a non-zero simplicial norm. It is obvious now that the image $h_n$ of $\bar{h}_n$ under the homomorphism induced by the universal central extension $\bar{A}_0 \longrightarrow A_0$ has a non-zero simplicial norm, and its image under the canonical isomorphism $H_n(BA_0) \longrightarrow H_n(BA_0^+)$ is spherical. $\square$

**Remark 1.** Theorem 6.6.D' in [Gr 2] easily implies that $A_0$ satsfies (5'.1). It is is not difficult to deduce (5'.2) from the fact that $A_0$ is torsion-free.

**Remark 2.** The groups $\tilde{A}_0$ were defined in [Cl] as follows. (We will describe the Clozel definition only in the situation $p = 2n - 1$, $q = 1$ which is of interest to us.) Let $N = p + q = 2n$. Let $F$ be a totally real number field and $F_c$ be its totally imaginary quadratic extension (e.g. $F = \mathbb{Q}(\sqrt{2})$ and $F_c = F(\sqrt{-3})$). Let $D$ be a division algebra of degree $N^2$ over $F_c$. Assume that $D$ satisfies the following condition (R): At any place $v$ of $F_c$, $D_v = D \bigotimes_{F_c} F_{c,v}$ is either isomorphic to $M_n(F_{c,v})$ or a division algebra. Assume further that $D$ is endowed with an involution of the second kind $x \longrightarrow x^*$, an involutive antiautomorphism acting on $F_c \subset D$ by conjugation. Consider the unitary group $U(F) = \{d \in D : dd^* = 1_N\}$. Assume further that (K): at all infinite places of $F$ except one $U(F)$ is compact and isomorphic to $U(n)$, and at the remaining infinite place $U(F)$ is isomorphic to $U(N-1, 1) = U(2n-1, 1)$. By restriction of scalars, the group $U/F$ defines a $\mathbb{Q}$-group denoted by $G_0$. (Thus, $G_0(\mathbb{Q}) = U(F)$.) Now $G_0(\mathbb{R}) = K \times U(2n-1, 1)$, where $K$ is a compact group. Let $\Gamma$ be a cocompact lattice in $U(2n-1, 1)$ commensurable with the projection of an arithmetic subgroup of $G_0(\mathbb{Q}) = U(F)$. Then Clozel proves that it is always possible to find a torsion-free congruence subgroup $\tilde{A}_0$ of $\Gamma$ with the required cohomology properties.

The existence for any $n \geq 2$ of a division algebra $D$ satisfying (R) endowed with an involution of the second kind such that the corresponding unitary group $U(F)$ satisfies the condition (K) is treated in [Cl] as a known fact but not explicitly stated. This fact is the immediate corollary of Proposition 2.3 in [Cl0]. The proof of Proposition 2.3 in [Cl0] is based on the results from [Kt0] and [Kt1]. (We are grateful to R. Kottwitz for explaining to one of us this proof in detail.)

**Remark 3.** Note that the second construction combined with the construction of the homology sphere $\Sigma^n$ at the beginning of Step 2 in section 2 imediately implies the following theorem:

**Theorem 10.** For any $n \geq 5$ there exists a smooth $n$-dimensional homology sphere of non-zero simplicial volume.

(Of course, for even $n$ Theorem 10 can be proved using the simpler first construction instead of the second construction.) Moreover, using the group $A_0$ obtained using the second construction in the proof of Theorem 1 in section 2 we obtain a (stronger) version of Theorem 1 where the property (2) is replaced by the following property (2'): "This manifold has simplicial volume not less than 1". (Of course, one can replace here 1 by any positive number.)

**Third construction.** Assume first that $n \geq 6$. (The case $n = 5$ will be considered at the end of the description of the construction). Let $K = \mathbb{Q}(\sqrt{2})$, $F = x_1^2 + \ldots + x_n^2 - \sqrt{2}x_{n+1}^2$ and $G = SO_F(n+1, K)^0$ be the group of unimodular $(n+1) \times (n+1)$ matrices with entries in the field $K$ preserving the form $F$ and mapping the upper half of the hyperboloid defined by the equation $F = -1$ into itself. Let $J$ be the ring of integers of $K$ and $L = J[\frac{1}{5}, \frac{1}{13}]$. Let $G(L) = G \bigcap GL_{n+1}(L)$. The group $G(L)$ is S-arithmetic and therefore finitely presentable. Consider the group $G_S$ defined as a subgroup of the adèle group of $G$ consisting of adèles whose $v$-component is equal to identity for all places $v$ but $5, 13$ and $\infty$. Thus, $G_S = G(\mathbb{R}) \times G(K_5) \times G(K_{13})$. Here $K_5$, $K_{13}$ denote as usual the result of the completion of $K$ with respect to the $p$-adic valuations of $K$ corresponding to $p = 5$ and $p = 13$, respectively. It is clear that $G(L)$ is a subgroup of $G_S$. A theorem of Borel and Harish-Chandra (cf. e.g. [Ma], Theorem I.3.2.7 and IX.1.4) implies that $G(L)$ is an irreducible cocompact lattice in $G_S$. Since 5 and 17 are primes of the form $4k+1$ the equation $x^2 + 1 = 0$ has solutions in $\mathbb{Q}_5$ and $\mathbb{Q}_{13}$. Therefore the ranks of $G(K_5)$ and $G(K_{13})$ are equal to $[(n+1)/2]$. Following [BorWal],XIII.2.1, let $r = r_\infty + r_f = 2[(n+1)/2] + 1$ be the rank of $G_S$. According to [BorWal], Theorem XIII.3.5 and Lemma XIII.3.2 (see also the proof of Proposition XIII.3.6 in [BorWal]) $H^q(G(L), \mathbb{R}) = H^q_{ct}(G(\mathbb{R}), \mathbb{R})$ if $q < r$. $H^*_{ct}(G(\mathbb{R}), \mathbb{R})$ denotes here the continuous cohomology of $G(\mathbb{R})$ which is known to be equal in our situation to the cohomology of the compact dual of the symmetric space $SO(n,1)^0/SO(n)$, that is of $SO(n+1)/SO(n) = S^n$. (Here $SO(n,1)^0$ denotes the connected component of $SO(n,1)$ containing the identity.) Thus, $H^q(G(L), \mathbb{R}) = 0$, if $q \in \{1, \ldots, n-1\}$ and $H^n(G(L), \mathbb{R}) = \mathbb{R}$. Moreover, these cohomology properties will hold for congruence subgroups of $G(L)$. According to a theorem of Selberg (cf. e.g. [Ser]) one can always find a torsion-free congruence subgroup of $G(L)$. Denote some such subgroup by $A_*$. Thus, $H_q(A_*, \mathbb{R}) = 0$, if $q \in \{1, \ldots, n-1\}$ and $\mathbb{R}$, if $q = n$. Denote by $\bar{h}_n \in H_n(BA_*)$ any nontrivial homology class such that its image in $H_n(BA_*, \mathbb{R})$ is a generator. As in the other constructions we can get rid of the first homology group of $A_*$ by forming an appropriate amalgamated free product with several copies of the acyclic group $< a, b, c, v_j | (9) - (12) >$, where the amalgamated subgroups are cyclic. Denote the resulting group by $A_0$ and the image of $\bar{h}_n$ in $H_n(\bar{A}_0)$ by $\tilde{h}_n$. The Serre generalization of Hurewicz theorem implies that for some non-zero integer $k$ and $h_n = k\tilde{h}_n$ the class $h_n^+ \in H_n(BA_0^+)$ is spherical. Condition (5'.2) of Proposition 4 is obvious, and it remains to check (5'.1). Lemmae 5.1-5.8 imply that it is sufficient to check (5'.1) for $(A_*, \bar{h}_n)$.

The group $G = SO_F(n+1, K)^0 \subset SO_F(n+1, \mathbb{R})^0$ is a group of (orientation preserving) isometries of the hyperbolic space $H^n$ which can be identified with the upper half of the hyperboloid $\{x \in \mathbb{R}^{n+1} | F(x) = -1\}$. Observe that since $A_* \subset G$ is a cocompact lattice in $G_S$, all its elements are semisimple (cf. [Ma],p.63). Therefore isometries of $H^n$ in $A_*$ cannot be parabolic (cf. [AVS] or [BGS]). Since $A_*$ is torsion-free they cannot be elliptic. Hence $A_*$ regarded as a group of isometries of $H^n$ conists only of hyperbolic isometries. Any almost nilpotent subgroup of $A_*$ consists of isometries preserving a specific geodesic $l$. Moreover, its normalizer in $A_*$ also consists only of isometries preserving the same geodesic *l*. The intersection of normalizers of almost nilpotent subgroups of $A_*$ preserving different geodesics is trivial. Any geodesic $l$ is determined by two "infinite" points $u$ and $v$. Denote by $A_{[u,v]}$ a subgroup of $A_*$ generated by all isometries preserving $[u,v]$. It is a subgroup of $G(L)_{[u,v]}$ formed by all elements of $G(L)$ preserving $[u,v]$ and of $SO_F(n+1,\mathbb{R})_{[u,v]}$ formed by all orientation preserving hyperbolic isometries preserving $[u,v]$. The last group is the product of $\mathbb{R}^\times$ and $SO(n-1)_{[u,v]}$, where $\mathbb{R}^\times$ is the multiplicative group of positive real numbers which can be identified with the group of dilations of $[u,v]$ and $SO(n-1)_{[u,v]}$ denotes the group of rotations around $[u,v]$. Denote the kernel of the homomorphism of $A_{[u,v]}$ to $\mathbb{R}^\times$ by $S_{[u,v]}$ and the image of this homomorphism by $R_{[u,v]}$. We see that $A_{[u,v]}$ isomorphic to the product of $S_{[u,v]}$ and the abelian group $R_{[u,v]}$. Almost nilpotent subgroups of $A_{[u,v]}$ are just products of almost nilpotent subgroups of $S_{[u,v]}$ with $R_{[u,v]}$. The normalizers of almost nilpotent subgroups of $A_{[u,v]}$ are products of normalizers of almost nilpotent subgroups of $S_{[u,v]}$ (in $S_{[u,v]}$) with $R_{[u,v]}$. Thus, $Q_{A_{[u,v]}}//A_{[u,v]}$ is the product of $Q_{S_{[u,v]}}//S_{[u,v]}$ and $BR_{[u,v]}$. Further, two almost nilpotent subgroups of $A_{[u_1,v_1]}$ and $A_{[u_2,v_2]}$ can be conjugate if and only if there is an isometry in $A_*$ which transforms the geodesic $[u_1,v_1]$ into $[u_2,v_2]$. In this case all almost nilpotent subgroups of $S_{[u_1,v_1]}$ are conjugate to the corresponding almost nilpotent subgroups of $S_{[u_2,v_2]}$. At this stage we see that $Q_{A_*}//A_*$ is made of disjoint copies of $Q_{S_{[u,v]}}//S_{[u,v]} \times BR_{[u,v]}$, where $[u,v]$ runs over a certain set of geodesics of $H^n$, all groups $R_{[u,v]}$ are abelian groups isomorphic to subgroups of $\mathbb{R}^\times$, and $Q_{S_{[u,v]}}//S_{[u,v]}$ is made of classifying spaces of subgroups of $S_{[u,v]}$. The image of the projection homomorphism $H_*(Q_{A_*}//A_*) \longrightarrow H_*(A_*)$ will be the union of homomorphisms corresponding to connected components of $Q_{A_*}//A_*$. Each of these last homomorphisms factors through $H_*(BS_{[u,v]} \times BR_{[u,v]})$ for the corresponding $[u,v]$. Thus, it is sufficient to show that $\bar{h}_n$ is not in the image of any inclusion homomorphism $\lambda : H_n(BS_{[u,v]} \times BR_{[u,v]}) \longrightarrow H_n(BA_*)$. We start from the observation that after tensoring with reals $\bar{h}_n$ becomes dual to a nonzero continuous cohomology class. Therefore, it remains nonzero under the inclusion homomorphisms $i^{ex} : H_n(A_*) \longrightarrow H_n(G(L^{ex}))$, where $L^{ex}$ denotes the ring of S-integers of any algebraic extension of $K = \mathbb{Q}(\sqrt{2})$ and the set of finite places $S$ is any finite set including 5 and 13. Choose $L^{ex}$ such that: 1) The group $SO(n-1, L^{ex})_{[u,v]}$ formed by all rotations around $[u,v]$ given by matrices with entries in $L^{ex}$ is isomorphic to $SO(n-1, L^{ex})$; and 2) The number of added finite places is large enough to ensure that $H_i(SO(n-1, L^{ex}), \mathbb{R}) = 0$ for all $i \in \{1,\ldots,n\}$ by virtue of Proposition XIII.3.6. (It is sufficient to add just two places 17 and 29 to finite places 5 and 13 which are already there.) Now note that if $\bar{h}_n$ is in the image of $\lambda$, then $i^{ex}(\bar{h}_n)$ is in the image of $i^{ex} \circ \lambda$. The last homomorphism factors through $H_n(BSO(n-1, L^{ex})_{[u,v]} \times BR_{[u,v]})$. The last group is isomorphic to $H_n(R_{[u,v]})$ by virtue of the choice of $L^{ex}$.

It remains to prove that for any geodesic $[u,v]$ the image of $H_n(R_{[u,v]}, \mathbb{R})$ in $H_n(G(L^{ex}), \mathbb{R})$ under the homomorphism induced by the inclusion is trivial. (This will be clearly sufficient to conclude that (5'.1) holds.)

At this point one can complete the proof in two different ways. We give the more elementary. (The other is based on a more systematic exploitation of the above connec-tion between group cohomology and continuous cohomology for $G(L^{ex})$, $R_{[u,v]}$ and other auxiliary groups.) We are going to show that $H_n(R_{[u,v]}) = 0$. Note that $u$ and $v$ are eigenvectors of any matrix from $SO_F(n+1, \mathbb{R})_{[u,v]}$. The restriction of the corresponding linear transformation on the linear space spanned by $u$ and $v$ is $diag\{\lambda, \frac{1}{\lambda}\}$ in the basis $\{u,v\}$ for some real $\lambda$. Therefore, if such matrix belongs to $G(L)_{[u,v]}$ then $\lambda$ must be a unit in the ring of S-integral elements of a quadratic extension of $\mathbb{Q}(\sqrt{2})$ where the set of added finite places is $\{5, 13\}$. Thus, $R_{[u,v]}$ is isomorphic to a multiplicative group which consists entirely of such numbers. An elementary ad hoc calculation shows that if $\lambda_1$, $\lambda_2$ and $\lambda_1\lambda_2$ are such numbers (from possibly different quadratic extensions of $\mathbb{Q}(\sqrt{2})$), and $\lambda_1$, $\lambda_2$ are not square roots of units in the ring of S-arithmetic integers of $\mathbb{Q}(\sqrt{2})$, ($S = \{5, 13\}$), then $\lambda_1$ and $\lambda_2$ belong to the same quadratic extension of $\mathbb{Q}(\sqrt{2})$. Now it is easy to see that any finitely generated multiplicative group which consists entirely of such numbers contains a subgroup of finite index which is a subgroup of the group of units of the ring of S-algebraic integers of one quadratic extension of $\mathbb{Q}(\sqrt{2})$. Now the generalization of the Dirichlet unit theorem for the rings of S-integral elements of algebraic extensions of $\mathbb{Q}$ implies that any such finitely generated multiplicative group of $\mathbb{R}^\times$ contains a subgroup of finite index isomorphic to a subgroup of $\mathbb{Z}^5$ (cf. [FT], Theorem 38 in Ch. IV). Therefore the $q$-th homology group with real coefficients of any such group is trivial for any $q \geq 6$. Therefore the $q$-th homology group of $R_{[u,v]}$ is trivial for any $q \geq 6$.

In the case when $n = 5$ we define $L$ as $J[\frac{1}{5}]$ and proceed as above. ($J$ is, as before the ring of integers of $\mathbb{Q}(\sqrt{2})$.) Then at the end of the proof we will see that $H_q(R_{[u,v]})$ is trivial for any $q \geq 5$. This is sufficient to prove (5'.1) almost as above. (There will be an easily resolvable minor difficulty in the proof of the fact that $H_i(SO(4, L^{ex}), \mathbb{R}) = 0$ for $i = 1, \ldots, 5$ and an appropriate $L^{ex}$ due to the fact that $SO(4)$ is the product of $SO(3)$ and $S^3$, and we do not have the irreducibility required in the cited results of [BorWal], XIII.3. But, of course, we can represent our lattice as the product of lattices, and apply Proposition XIII.6 from [BorWal] to these lattices separately. To ensure the vanishing of the first five homology groups of $SO(4, L^{ex})$ with real coefficients we can add to $L^{ex}$ five finite places $p_i$, where $p_i$ are primes of the form $4k+1, k > 1$.)

However the first part of the proof will imply only the vanishing of the first three homology groups of $A_*$ with real coefficients (instead of four). But a generalization due to Serre of the Hurewicz theorem theorem implies that if $X$ is a simply connected CW-complex such that $H_2(X, \mathbb{Q}) = H_3(X, \mathbb{Q}) = 0$ then the Hurewicz homomorphism $\pi_5(X) \otimes \mathbb{Q} \longrightarrow H_5(X, \mathbb{Q})$ is surjective. This fact is sufficient to conclude that if $\tilde{h}_n \in H_n(A_*)$ is a non-trivial class dual to the cohomology class coming from the continuous cohomology group of $G(\mathbb{R})$, then for some its non-zero multiple $h_n$ the class $h_n^+ \in H_n(BA_*^+)$ is spherical.

**Remark.** Even in the case $n \geq 6$ we could define $L$ as $J[\frac{1}{5}]$ (instead of $J[\frac{1}{5}, \frac{1}{13}]$). In this case we will obtain the vanishing of homology groups of $A_*$ with real coefficients in all positive dimensions $\leq [(n+1)/2]$. But this is still sufficient in order to prove that $h_n^+ \in H_n(BA_*^+)$ is spherical, where $h_n$ is a nontrivial multiple of the homology class dual to the generator of $H_{ct}^n(A_*, \mathbb{R})$. (One can use rational homotopy theory (cf. e.g. [S], [GM] or [FHT]).)

## 4. Proof of Theorem A from Theorem 1.

As the reader might have guessed from the title of this section we will prove the following:

**Proposition 0.2. A.** If Theorem A, parts (i)-(iv) is false for some $n \geq 5$, then for the same value of $n$ Theorem 1.A and its generalization in part C for an arbitrary compact smooth $n$-dimensional manifold $M_0^n$ are false;

**B.** If Theorem A, parts (i)-(v), is false for some $n \geq 5$, then for the same value of $n$ Theorem 1.B and its generalization in part C of Theorem 1 for an arbitrary smooth compact $n$-dimensional manifold $M_0^n$ are false.

To prove Proposition 0.2 we are going first to describe how one can effectively construct for any fixed $n$ and any given positive rational $x$ and $v$ a finite set $Net(n,x,v)$ of $n$-dimensional Riemannian manifolds with the following properties:
(i) Any compact Riemannian manifold $W^n$ such that $|K| \leq 1$, its diameter does not exceed $x$ and its volume is at least $v$ is $\delta$-close (in the Gromov-Hausdorff metric) to an element $E$ of $Net(n,x,v)$, where $\delta$ is defined as follows. Let $conv(n,x,v)$ denote the lower bound for the convexity radius of manifolds such that $|K| \leq 3$, $vol \geq v/3$ and $diam \leq x+2$. implied by the Cheeger bound for the injectivity radius given in [C] (cf. [Ch], Theorem 7.6 and 7.9). For large $x$ $conv(n,x,v)$ behaves as $const(n)v/\exp((n-1)x)$. Let $\delta = conv(n,x,v)/(1000n^2)$. This choice of $\delta$ ensures that any two $10\delta$-close Riemannian manifolds such that $|K| \leq 3, vol > v/3$ and $diam \leq x+2$ are homotopy equivalent (cf. [P]). Also, this element $E$ of $Net(n,x,v)$ must have volume not exceeding $3\times$ the volume of the Riemannian manifold $W^n$.
(ii) Any element of $Net(n,x,v)$ is a Riemannian manifold such that the absolute value of its sectional curvature does not exceed 3, its volume is not less than $v/3$ and its diamter is not greater than $x+2$.
Of course, any element of $Net(n,x,v)$ must be representable in a finite form. (Otherwise it does not make sense to say that $Net(n,x,v)$ can be effectively constructed.) Here is an outline of the construction: First, we are going to show that any Riemannian metric $(M^n,g_0)$ on a compact $n$-dimensional manifold such that $|K_{g_0}| \leq 1$, $diam_{g_0}(M^n) \leq x$ and $vol_{g_0}(M^n) \geq v$ can be approximated in the $C^2$-norm by a "nice" Riemannian metric $g_{nice}$ on $M^n$ with the following properties: (a) $M^n$ is endowed with a structure of a $C^2$-smooth semialgebraic manifold "made of" coordinate charts such that their number does not exceed the value of a certain computable function of $x$ and $v$ (this function can be written down explicitly), and all transition functions are $C^2$-smooth semialgebraic functions of local coordinates. (Here and below "a computable function" means a function which either increases or decreases with any of its arguments and its restriction on the set of rational values of arguments is Turing computable. The dimension $n$ is always regarded as a constant.) There exist computable upper bounds (in terms of $x$ and $v$) for the minimal complexity of diagrams of the graphs of transition functions, and the absolute values of coefficients of polynomials entering some diagrams of the minimal complexity of graphs of transition functions. (Recall that a semialgebraic function is, by definition, a function such that its graph is a semialgebraic set. Any semialgebraic set has a (non-unique) diagram of the form $\bigcup_i \bigcap_j \{P_{ij} ?_{ij} 0\}$, where $P_{ij}$ are polynomials and for any $i, j$ $?_{ij}$ denotes one of the signs $>, <, =$. The complexity of the diagram is, by definition, $\Sigma_{i,j} \deg P_{ij}$.) Furthermore, any chart $B_R(0) \longrightarrow M^n$ can be extended to completely cover all overlapping images of local charts; the transition functions defined on the extended local charts still will be $C^2$-smooth semialgebraic functions satsfying the same computable upper bounds for the minimal complexity of diagrams and the size of coefficients of the polynomials in a diagram of the minimal complexity;

(b) $|K_{g_{nice}}| \leq 1.5$, $\frac{vol_{g_{nice}}(M^n)}{vol_{g_0}(M^n)} \in [0.5, 2]$, $diam_{g_{nice}}(M^n) \leq x+1$, $d_{GH}((M^n, g_0), (M^n, g_{nice})) \leq \delta$;

(c) The metric $g_{nice}$ is $C^2$-smooth semialgebraic in local coordinates. There exist computable upper bounds for the minimal complexity of diagrams for graphs of $g_{nice}$ on local charts as well as computable upper bounds for the absolute values of coefficients of polynomials entering some diagrams of the minimal complexity.

When we already know that such $g_{nice}$ exists, the problem reduces to constructing a required net in the space of such "nice" metrics on all compact $n$-dimensional manifolds. The last problem is clearly algorithmically solvable because of the following reasons: The set of "nice" metrics can be identified with a subset of a compact semialgebraic set of an Euclidean space. One is able to evaluate volume, diameter and $\sup |K|$ within to any accuracy, and all these quantitites are Lipschitz on this set with Lipschitz constants which can be effectively majorized. Moreover, we are allowed to include in the net Riemannian metrics satisfying slightly worse bounds for the sectional curvature, volume and diameter, which enables us to perform approximate calculations.

Thus, it is sufficient to describe how to find a "nice" Riemannian metric $g_{nice}$ on $M^n$ with the properties (a)-(c). As it was already mentioned, the classical Cheeger inequality implies a computable (in terms of $n$, $x$, $v$) lower bound for the injectivity and, as a corollary, convexity radius of any compact Riemannian manifold with $|K| \leq 3$, $vol \geq v/3$, $diam \leq x+2$ (cf. [Ch], Theorems 7.6 and 7.9). Denote this bound by $conv(n, x, v)$. Rescale the Riemannian metric $g_0$ to make all distances increase by the factor $1/conv(n, x, v)$. (Without any loss of generality we can assume that this number is greater than one. This rescaling and the inverse rescaling at the end of construction below are only the matter of convenience and not really necessary.) The absolute values of sectional curvatures will be bounded now by $conv(n, x, v)^2$, and the convexity radius will be $\geq 1$. Let us apply now the Ricci flow, ($\frac{dg}{dt} = -2Ric(g)$), to smooth out the resulting Riemannian metric. The classical result of Bemelmans, Min-Oo and Ruh ([BMR], [Ba], [F]) guarantees the existence of the solution for times not exceeding a certain positive constant $C_n$. For all $t$ not exceeding a certain positive constant $\bar{C}_n \leq C_n$ the absolute values of sectional curvature will not exceed 1 and the volume of any ball will not exceed two times its value for $t = 0$. Now the local injectivity radius estimate in [CGT] (Theorem 4.7) implies that for all values of $t$ up to a certain positive constant $\hat{C}_n \leq C_n$ the injectivity radius of $(M^n, g_t)$ will be not less than $\frac{1}{2}$. Choose a system of local harmonic coordinates on $(M^n, g_t)$ constructed as in [JK]. (One can define harmonic coordinates on any metric ball of radius $< \rho = const > 0$, where $const$ denotes here and below a (every time different) constant which does not depend on $x, v, g_0$ ([JK], Satz 5.1). Below we will be more specific about how exactly we choose the system of local harmonic coordinates. But at the moment the details do not matter.) For $t$ not exceeding some $T_0 = const > 0$ the supremum of the $C^0$-norm of the curvature tensor during this process will not exceed $100/99\times$ the supremum of the $C^0$-norm of the curvature tensor of $g_0$. Thus, absolute values of sectional curvature of metrics $g(t)$ for $t \leq T_0$ will be bounded from above by $1.4\, conv(n,x,v)^2$. Furthermore, $|g(t) - g_0|_{C^0} \leq const\, t$ (cf. [BMR]). Therefore choosing $t = T = const\, \min\{1, \frac{conv(n,x,v)}{x}\}$ (for an appropriate $const > 0$) ensures that $d_{GH}((M^n, g_0),(M^n, g(T))) \leq 1/(4000n^2)$, and that $vol_{g(T)}(M^n) \geq 3/4 vol_{g_0}(M^n)$. For $t \in (0,T_0]$ one has upper bounds for $|\nabla^k Rm(g(t))|$ of the form $const/t^{k/2}$ for $k = 1,2,3$ (as well as for any other fixed $k$), where $Rm(g(t))$ denotes the curvature tensor of $g(t)$ (cf. [Ba]). As a corollary, we get upper bounds for the $C^3$-norm of $g(T)_{ij} - \delta_{ij}$ on any coordinate chart, where $g_{ij}$ are regarded as functions of local harmonic coordinates ([BMR], eq. (7) on p. 72). Proceeding as in the proof of Lemma 4.3 in [C] we can obtain upper bounds for the $C^3$-norms of the transition functions between the local charts. Note that the inequalities (5.2) and (5.5) on p. 65 of [JK] imply that for some positive $\rho_1 = const < \rho$ and any point $p \in (M^n, g(T))$ the harmonic coordinates $H : B_{\rho_1}(p) \longrightarrow \mathbb{R}^n$ chosen in the metric ball of radius $\rho_1$ around $p$ have the following property: For any $r \in (\rho_1/100, \rho_1)$ $B_{9r/10}(0) \subset H(B_r(p)) \subset B_{10r/9}(0) \subset \mathbb{R}^n$. (The inequality (5.5) in [JK] expresses the fact that the differential of $H$ is sufficiently close to an isometry.) Now let us be more specific about the choice of a system of local harmonic coordinates. Choose $\rho_1/20$-net in $(M^n, g(T))$ such that the balls of radius $\rho_1/50$ centered at the points of the net are disjoint. Choose the harmonic coordinates $H$ in the metric ball of radius $\rho_1$ around every point of the net. Consider the restrictions of these local coordinates on (non-metric) balls $H^{-1}(B_{\rho_1/9}(0))$ of radius $\rho_1/9$ in local coordinates as a system of coordinate charts. For any two overlapping coordinate charts the local coordinates defined on one of the coordinate charts extend to the other. Also, the number of coordinate charts can be easily majorized using the Bishop-Gromov inequality. Let us now approximate in the $C^2$-norm all transition functions by polynomials (or, correspondingly, smooth algebraic functions inverse to the approximating polynomials) on balls in harmonic coordinates of radius $\rho_1/3$ centered at the origins of the harmonic coordinate charts. The availability of the upper bounds for the $C^3$-norms of the transition functions enables us to obtain an explicit upper bound for the degrees and the absolute values of coefficients coefficients of approximating polynomials in terms of the accuracy of the approximation, $x$ and $v$. (We will explain below how to choose the accuracy of the approximation.) Also, using any standard constructive proof of the Weierstrass approximation theorem we can get upper bounds for $C^3$-norms of the approximations of transition functions in terms of $x$ and $v$ but not the accuracy of the approximation. At this stage our manifold acquires a structure of a $C^2$-smooth semialgebraic manifold. Observe, however that the expressions for the Riemannian metric $g(T)$ in the local coordinates on the overlaps of coordinate charts will not be quite the same when it is recalculated from a different coordinate chart using new transition functions. Still, the $C^2$-norm of the difference can be majorized in terms of the $C^2$-norm of the error of approximation of transition functions. Now let us approximate the expressions for the Riemannian metric in local coordinates in the balls of radius $\rho_1/3$ around the origin by polynomials in the $C^2$-norm. Since we have the control over $C^3$-norms we can effectively majorize the degrees and the absolute values of the approximating polynomials in terms of the accuracy of approximation which will be defined below, $x$ and $v$ (recall that $n$ is regarded as a constant). We approximate $g_{ij}$ and $g_{ji}$ by the same polynomial. Choose a $C^2$-smooth semialgebraic partition of unity corresponding to the chosen system of coordinate charts, such that the $C^2$-norms of all partition functions, complexity of diagrams of the graphs of all partition functions as well as the absolute values of coefficients of polynomials in a diagram of minimal complexity (for any local system of coordinates) do not exceed a certain computable bound depending on $x$ and $v$ (cf. [BCR]). We demand also that these partition functions must vanish outside of balls of radius $\rho_1/9$ (in harmonic coordinates) around the centers of coordinate charts. (These conditions clearly can be satisfied since already balls of radius $\rho_1/18$ around the centers of coordinate charts cover all manifold. These partition functions can be defined, for example, using the following explicit formulae. Let $Ch_1, \ldots, Ch_N$ be coordinate charts and $h_1^{(k)}, \ldots, h_n^{(k)}$ be local coordinates in $Ch_k$. Let $\psi_k(x) = (\sum_{i=1}^n (h_i^{(k)}(x))^2 - (\frac{\rho_1}{9})^2)^3$, if $x \in Ch_k$ and $\sum_{i=1}^n (h_i^{(k)}(x))^2 \leq (\frac{\rho_1}{9})^2$, and $\psi_k(x) = 0$, otherwise. Define $\lambda_k(x)$ as the ratio $\frac{\psi_k(x)}{\sum_{i=1}^N \psi_i(x)}$. The functions $\lambda_1, \ldots, \lambda_N(x)$ constitute the required $C^2$-smooth semialgebraic partition of unity.) Using this partition of unity we can glue the algebraic approximations to the Riemannian metric $g(T)$ into a $C^2$-smooth semialgebraic metric on the manifold. That is, this will be the case if we choose the accuracies of the approximation of transition functions and of the Riemannian metric (on the previous steps) to ensure the positive definiteness of the metric. Also, we want to ensure that the Gromov-Hausdorff distances between the resulting Riemannian manifold and the Riemannian manifold $(M^n, g(T))$ must not exceed $1/(4000n^2)$, and that the volume is changed not more than by the factor of $3/4$. Clearly this can be done. Now after the rescaling back by the factor $conv(n, x, v)$ we will get the desired Riemannian metric $g_{nice}$.

Suppose now that we are given a smooth $n$-dimensional homology sphere (or, corespondingly, homology $M_0^n$) which is either diffeomorphic to $S^n$ or satisfies properties (1)-(3) from the text of Theorem 1 (correspondingly (1), (2) if $M_0^n$ is spin or just (2) if $M_0^n$ is not spin). Without any loss of generality we can assume that the given manifold is endowed with a Riemannian metric. Assume now that $S^n$ (or, more generally, the manifold $M_0^n$ from the text of Theorem 1, part C) admits Riemannian metrics of arbitrarily small volume such that $|K| \leq 1$. (For $S^n$ this condition is equivalent to the condition that $n$ is odd.) Using the nets $Net(n, x, v)$ we can design an algorithm based on the assumption that Theorem A, part (i)-(iv), is false (for $M^n = S^n$ or correspondingly $M^n = M_0^n$) and bringing to contradiction Theorem 1.A and its generalization (in part C) for any such manifold $M_0^n$. Take $v(n) = w(n)/100$, where $w(n)$ is the constant from the text of Theorem 1, and $v(n)$ is the constant from the text of Theorem A. Note that the assumption that Theorem A is false implies that for any effectively majorizable function $\phi$ and for all sufficiently large $x$ any Riemannian metric on $S^n$ (correspondingly, $M_0^n$) such that $|K| \leq 1$ and $diam \leq x$ can be connected by a sequence of not very large "jumps" (that is, "jumps" less than $\delta$ defined in the part (i) of the definition of $Net(n, x, v)$ but for $Net(n, \phi(x), v(n))$) via Riemannian metrics satisfying $|K| \leq 1$ and of diameter $\leq \phi(x)$ with a Riemannian metric of volume $\leq v(n)$. We are going to proceed as follows: Find an upper bound $d$ for the diameter and an estimate for the volume of the given Riemannian manifold such that the relative error does not exceed 0.1. Let $v$ be any fixed positive rational number between $w(n)/3$ and $w(n)/2$. If the found estimate for the volume of the manifold is less than $v$ then the manifold is $M_0^n$. Compute an upper bound $X$ for $\phi(d)$. Construct $Net(n, X, v/100)$. We want to check whether or not there exist a finite sequence of "jumps" of length $\leq 3\delta = 3conv(n, X, v/100)/(1000n^2)$ via elements of the net $Net(n, X, v/100)$ connecting the considered Riemannian structure on the manifold and a Riemannian structure on the manifold of volume less than $v/9$. The manifold is $S^n$ (or, correspondingly, $M_0^n$) if and only if such a sequence exists. Moreover, this statement will remain true if we replace $3\delta$ by $4\delta$ and $v/9$ by $v/5$. This gives us some room for an error in approximate calculations. Also note that we can replace the original Riemannian manifold by (any) its approximation in $Net(n, X, v/100)$ constructed as above, and all these observations still remain true. But now the algorithmic problem becomes trivial: It is sufficient first to represent $Net(n, X, v/100)$ as a graph such that its set of vertices coincides with the set of elements of $Net(n, X, v/100)$ and two vertices are connected by an edge if and only if the Gromov-Hausdorff distance between corresponding metric spaces computed using any fixed approximated algorithm within to the accuracy $\delta/3$ does not exceed $3.5\delta$, and then to check whether or not the approximation to the given Riemannian manifold is in the same component of this graph as a vertex corresponding to a Riemannian manifold such that the approximate value of its volume calculated using some fixed algorithm within to the accuracy $v/100$ does not exceed $v/7$. Denote this graph by $Gr(n, X, v/100)$. It will be used also in the proof of part B of Proposition 0.2.

The case, when any Riemannian metric on $S^n$ (or $M_0^n$) such that $|K| \leq 1$ has a volume greater than or equal to $v_0$ for some positive $v_0$ is much easier and can be treated as above if $v_0 < v/30$ or as in [N3] if $v_0 \geq v/30$. Indeed, in the last case assume first that $M_0^n$ is simply-connected. The assumption that Theorem A, parts (i)-(iv), is false is equivalent to the assumption that for all sufficiently large $x$ any two Riemannian structures on $S^n$ (corresp. on $M_0^n$) of $|K| \leq 1$ and of diameter $\leq x$ can be connected by a sequence of short "jumps" in $R_1(S^n)$ (or $R_1(M_0^n)$) with a "controllable" increase of diameter. Then we can compute as in [N3] a number $Fl$ such that any loop of length $\leq l$ can be contracted to a point via loops of length $\leq Fl \times l$. The knowledge of such $Fl$ will then enable us to check whether or not the given homology $M_0^n$ for which it is known that it is either diffeomorphic to $M_0^n$ or nonsimply-connected is simply-connected. (see [N3] for the details). To complete the proof of Proposition 0.2.A we must also consider the case when $M_0^n$ is nonsimply-connected and $|K| \leq 1$ implies $volume \geq v_0$ for some positive $v_0 \geq v/30$. This case can be treated as the case of a simply-connected $M_0^n$ such that $|K| \leq 1$ implies $volume \geq v_0 \geq v/30$ but with the following modification: We have a sequence of manifolds with fundamental groups either isomorphic to $\pi_1(M_0^n)$ (in which case the manifold is diffeomorphic to $M_0^n$) or is the free product of $\pi_1(M_0^n)$ and a non-trivial group. We leave the proof of this case following the lines of [N3] as an exercise to the reader, but note that it formally follows from the proof of the part B of Lemma 0.2 given below.

In order to prove part B of Proposition 0.2 assume that Theorem A, parts (i)-(v), is false for some manifold $M_0^n$ of dimension $\geq 5$. Consider the sequence $\{M_i^n\}_{i=1}^\infty$ of effectively constructed in the proof of Theorem 1 smooth Riemannian manifolds with the same homology groups as $M_0^n$. ($M_i^n$ is diffeomorphic to the connected sum of $M_0^n$ and $S_i$ constructed on Steps 2,3 in Section 2.) Our goal is to exhibit an algorithm solving the algorithmic problem of recognition of $M_0^n$ among the manifolds $M_0^n$ in a recursively growing with $N$ time and using only a sublinearly growing with $N$ number of bits of oracle information to solve the recognition problem for the first $N$ elements of the sequence of manifolds ($N$ is, of course, variable). First, describe the oracle information which our algorithm will be using. We request a description of a representative from every "jump" component of $diam^{-1}((0, const_1(n)(\ln N)^{\frac{1}{n}}])$ (where $diam$ is regarded as a functional on $R_1(M_0^n)$) such that the volume is bounded from below by $v/30$ on this "jump" component, and which contains at least one manifold of convexity radius $\geq 0.5$ and volume $\leq const_3(n)\ln N$. This representative must have convexity radius $\geq 0.5$ and volume $\leq const_3(n)\ln N$. Here $const_1(n), const_3(n)$ are the constants from the text of Theorem 1, part B. By a "jump" component we mean the equivalence class of Riemannian structures from $diam^{-1}((0, const_1(n)(\ln N)^{\frac{1}{n}}]$ with respect to the equivalence relation obtained as the transitive closure of the relation "to be $E(n)\exp(-(n-1)const_1(n)(\ln N)^{\frac{1}{n}})$-close in the Gromov-Hausdorff metric", where $E(n)$ is the positive constant from the text of Theorem A. Denote by $j(N)$ the number of such "jump" components. Obviousy, $j(N)$ is a non-negative integer number not exceeding $J$ (defined in the text of Theorem A) for $x = const(n)(\ln N)^{\frac{1}{n}}$. (Here we assume that the partition $R_i(M_0^n, x)$, $i = 1, \ldots, I$, described in the text of Theorem A is the finest possible: Any two Riemannian structures in $R_i(M_0^n, x)$ belong to the same "jump" component of $diam^{-1}((0, \phi(x)])$. The finiteness of $j(N)$ follows from the precompactness of the subset of $R_1(M_0^n)$ formed by Riemannian structures having any common upper bound for diameter.) Since Theorem A, (i)-(v) is false by our assumption, then for any positive $c$ there exist arbitrarily large $x$ such that $J \leq \exp(cx^n) \leq N^{c\ const(n)}$. In particular, there exists an unbounded increasing sequence of values of $x$ and, thus, of $N$ such that $J \leq \sqrt{N}$. Let $\Delta = conv(n, x, v/100)/(1000n^2) \sim const(n)\exp(-(n-1)x)$ as $x \longrightarrow \infty$ while $v$ is fixed. We are going to present any chosen representative from any considered "jump" component by its $\Delta/10$-approximation by a finite metric space defined as follows: It is a $\Delta/10$-net in the Riemannian manifold such that balls of radius $\Delta/40$ around points of the net do not intersect. The number $L$ of points of the net is $\sim const(n)\frac{x^n}{\delta^n} \sim const(n)x^n \exp(n(n-1)x)$, and the number of distances between these points is $\sim const(n)x^{2n} \exp(2n(n-1)x)$. We would like to represent every length with the accuracy $\Delta/(1000L)$ This will require $\sim \ln \frac{L}{\Delta} \sim const(n)x$ binary digits. Summarizing, we need $\leq Const(n)x^{2n+1} \exp(2n(n-1)x) \leq Const_*(n) \exp(Const_{**}(n)(\ln N)^{\frac{1}{n}}))$ bits to present every representative (for some positive constants $Const(n), Const_*(n), Const_{**}(n)$). Thus, for any positive constant $C$ for infinitely many values of $N$ the total amount of the oracle information we are going to use does not exceed $N/C$. Thus, if we will be able to exhibit an algorithm using this oracle information and solving the recognition problem in time bounded by a computable function of $N$ we will obtain the desired contradiction.

For any of the given by the oracle metric spaces we can find in a recursive time all elements of the $Net(n, x, v/100)$ which are $conv(n, x, v/100)/(1000n^2)$-close to this metric space. Denote these elements by $E_i$. Then the algorithm performs a calculation of an upper bound for the diameter $d$ of the given manifold and then of an upper bound $X$ for $\phi(d)$. Find $Net(n, X, v/100)$ proceeding as described in the proof of part A of Proposition 0.2 above. Find any elements of $Net(n, X, v/100)$ which is $conv(n, x, v/100)/(500n^2)$-close to the given Riemannian manifold. Consider it as a vertex $A$ of the graph $Gr(n, X, v/100)$ (see the proof of part A of the proposition above for the definition of $Gr(n, X, v/100)$). The next step of the algorithm is the determination whether or not this vertex $A$ is in the same component of $Gr(n, X, v/100)$ as one of the vertices corresponding to the elements $E_i$ in $Net(n, x, v/100)$. (Recall that these elements are close to the finite metric spaces provided by the oracle). Of course, if the answer is "yes" then the manifold is $M_0^n$. If the manifold $M_0^n$ is such that any Riemannian metric on $M_0^n$ such that $Ric \geq -(n-1)$ has volume not less than $v/30$, then the answer "no" means that the given Riemannian manifold is not diffeomorphic to $M_0^n$. If $M_0^n$ does not have this property, and the answer is "no", then the given Riemannian manifold is diffeomorphic to $M_0^n$ if and only if the vertex $A$ is in the same component of $Gr(n, X, v/100)$ as a vertex of $Gr(n, X, v/100)$ corresponding to a Riemannian manifold of volume $\leq v/9$ if and only if the vertex $A$ is in the same component of $Gr(n, X, v/100)$ as a vertex corresponding to a Riemannian manifold of volume $\leq v/5$. Thus, it is sufficient to check whether or not $A$ is in the same component as a vertex corresponding to a Riemannian manifold such that the calculated value of its volume is $\leq v/7$, where the accuracy of the calculations of volume is $v/100$. Now it is obvious that there exists an algorithm checking this property (exactly as at the end of the proof of the part A of the proposition). $\square$

**Proof of Remark 2 after the text of Theorem A.** (1) To prove the strengthening of Theorem A stated as the first part of Remark 2 after the text of Theorem A we would like to prove a modified version of Theorem 1.B, where the condition that $Ric_{(S_i,\mu)} \geq -(n-1)$ implies $vol_\mu(S_i) \geq v(n)$ for any Riemannian metric $\mu$ on any of those homology spheres $S_i$ which are nonsimply-connected is replaced by the stronger condition that the constraint $Ric_{(S_i,\mu)} \geq -(n-1)$ for any Riemannian metric $\mu$ on nonsimply-connected $S_i$ implies $vol_\mu(S_i) \geq \exp(s(n)(\ln(i+1))^{\frac{1}{n}})$, where $s(n)$ is a positive constant but the condition (3) from the text of Theorem 1.A is dropped. (The proof of the corresponding version of Proposition 0.2 will be absolutely analogous to the proof of Proposition 0.2 in section 4.) To obtain this strengthened version of Theorem 1.B we can proceed exactly as above providing that we will be able to construct for any $j$ a smooth homology sphere $\Sigma^n(j)$ such that it admits a Riemannian metric with $\sup|K| \sim 1$ and of diameter between $c_1(n)j$ and $c_2(n)j$, of volume $\leq \exp(c_4(n)j)$ and such that its simplicial volume is not less than $\exp(c_3(n)j)$ for some positive constants $c_1(n)$, $c_2(n)$, $c_3(n)$ and $c_4(n)$. These homology spheres $\Sigma^n(j)$ for $j = [(\ln(i+1))^{\frac{1}{n}}]$ will then be used in the construction of $S_i$ on Step 2 of the proof of Theorem 1 in Section 2 instead of the homology sphere $\Sigma^n$ defined there. Observe, that in both the original and the modified constructions we will obtain $(S_i, \mu_i)$ such that $diam_{\mu_i}(S_i) \sim const(n)(\ln(i+1))^{\frac{1}{n}}$. However, in the modified construction the volumes $vol_{\mu_i}(S_i) \sim \exp(const(n)(\ln(i+1))^{\frac{1}{n}})$, which is worse than the upper bound for the volume $const_3(n)\ln(i+1)$ in the original construction. Still, the examination of the proof of Theorem A in section 4 shows that this weaker upper bound for the volumes of $(S_i, \mu_i)$ in Theorem 1.B is sufficient. Assume that a group $A_0$ has the properties (1)-(5) stated in the text of Proposition 4 at the beginning of section 2. (In the previous section we proved that for any $n$ such groups exist.) Consider the homology sphere $\Sigma^n$ such that its fundamental group is the universal central extension of $A_0$ introduced at the beginning of Step 2 in Section 2. If $A_0$ satisfies (5) then the simplicial norm of $\Sigma^n$ is non zero. Fix any Riemannian metric on $\Sigma^n$. Let $\Sigma^n(1) = \Sigma^n$, $\Sigma^n(2)$ be obtained from $\Sigma^n$ by taking the Riemannian direct sum with two copies of $\Sigma^n$, and more generally for any $j$ $\Sigma^n(j)$ is obtained by taking a Riemannian connected sum of $1+2+\ldots+2^{j-1}$ copies of $\Sigma^n$ along the binary tree of height $j$. It is clear that the diameter of $\Sigma^n(j)$ with such Riemannian metrics will grow linearly with $j$, and the simplicial volume of $\Sigma^n(j)$ will be equal to $(2^j-1)\times$ the simplicial volume of $\Sigma^n$.;

(2) The strengthening of Theorem A stated as the second part of Remark 2 after the text of Theorem A immediately follows from our proof of Theorem A and the Injectivity Radius Estimate at the bottom of p. 14 of [Gr 1]. (We would like to thank Misha Gromov who attracted our attention to this injectivity radius estimate and the possibility of its application in the context of our work.)

**Acknowledgements.** We are grateful to J. Arthur, S. Cappell, J. Cheeger, M. Gromov, B. Khesin and R. Kottwitz for helpful discussions.

### References.

A.N.: Department of Mathematics, University of Toronto, Toronto, Ontario, M5S 3G3, Canada; e-mail address: alex@math.toronto.edu

S.W.: Department of Mathematics, University of Chicago, Chicago, Illinois, 60637-1538, USA; e-mail address: shmuel@math.uchicago.edu